\documentclass[12pt,reqno]{amsart}
\usepackage{amsthm}
\usepackage{amsmath}
\usepackage{version}
\usepackage{graphicx}
\usepackage{enumerate}
\usepackage[usenames,dvipsnames] {color}
\numberwithin{equation}{section}

\def\d{\mathbb{D}}
\def\c{\mathbb{C}}
\def\t{\mathbb{T}}
\def\r{\mathbb{R}}

\newcommand\calc{\mathcal{C}}

\def\calr{\mathcal{R}}
\def\cale{\mathcal{E}}
\def\cala{\mathcal{A}}
\def\calf{\mathcal{F}}

\def\calv{\mathcal{V}}
\def\calw{\mathcal{W}}

\def\der#1#2{\frac{\partial#1}{\partial#2}}

\def\be{\begin{equation}}
\def\ee{\end{equation}}

\def\s0{s_0}
\def\p0{p_0}
\newcommand\ev{\mathrm {ev}}

\newcommand\calo{\mathcal{O}}

\DeclareMathOperator{\Kob}{{\mathrm Kob}}
\DeclareMathOperator{\orb}{{\rm Orb}}

\DeclareMathOperator\arctanh{arctanh}

 \newtheorem{theorem}{Theorem}[section]
 \newtheorem{corollary}[theorem]{Corollary}
 \newtheorem{lemma}[theorem]{Lemma}
 \newtheorem{proposition}[theorem]{Proposition}

\newtheorem{definition}[theorem]{Definition}
\newtheorem{remark}[theorem]{Remark}

\newtheorem{fact*}{Fact}

\DeclareMathOperator{\diag}{diag}
\DeclareMathOperator\re{Re}
\DeclareMathOperator\im{Im}
\DeclareMathOperator\spa{span}
\DeclareMathOperator\rank{rank}

\DeclareMathOperator\aut{Aut}
\DeclareMathOperator\laut{\mathrm{Lie}(\mathrm{Aut}~\mathbb{D})}
\newcommand\half{{\tfrac 12}}

\newcommand\dd{\mathrm d}
\newcommand\id{{\mathrm{id}}}
\newcommand\idd{\mathrm{id}_\mathbb{D}}
\newcommand\e{\mathrm e}
\newcommand{\eps}{\varepsilon}

\newcommand{\T}{\mathbb{T}}

\newcommand{\D}{\mathbb{D}}
\newcommand{\C}{\mathbb{C}}

\newcommand{\set}[1]{\left\{#1\right\}}

\newcommand{\ran}[1]{\operatorname{ran}#1}

\newcommand{\ip}[2]{\left\langle #1, #2 \right\rangle}

\newcommand{\inv}{^{-1}}

\newcommand{\ph}{\varphi}
\renewcommand\phi{\varphi}
\newcommand\al{\alpha}
\newcommand\ga{\gamma}
\newcommand\Ga{\Gamma}
\newcommand\de{\delta}

\newcommand\la{\lambda}
\newcommand\La{\Lambda}
\newcommand\up{\upsilon}
\newcommand\ups{\upsilon}

\newcommand\beq{\begin{equation}}

\newcommand\eeq{\end{equation}}
\newcommand\df{\stackrel{\rm def}{=}}

\newcommand\black{\color{black}}
\newcommand\red{\color{red}}

\newcommand\nn{\nonumber}
\newcommand\bbm{\begin{bmatrix}}
\newcommand\ebm{\end{bmatrix}}
\newcommand\bpm{\begin{pmatrix}}
\newcommand\epm{\end{pmatrix}}
\numberwithin{equation}{section}
\newcommand\nin{\noindent}
\makeindex

\newtheorem{lem}[theorem]{Lemma}
\newtheorem{prop}[theorem]{Proposition}
\newtheorem{thm}[theorem]{Theorem}

\newtheorem{defin}[theorem]{Definition}

\begin{document}
\title[A Characterization of the Symmetrized Bidisc]{A Geometric Characterization of the Symmetrized Bidisc}
\author{Jim Agler}
\address{Department of Mathematics, University of California at San Diego, CA \textup{92103}, USA}
\thanks{Partially supported by National Science Foundation Grants
DMS 1361720 and 1665260, a Newcastle URC Visiting Professorship and the Engineering and Physical Sciences Research Council grant EP/N03242X/1 }

\author{Zinaida Lykova}
\address{School of Mathematics, Statistics and Physics, Newcastle University, Newcastle upon Tyne
 NE\textup{1} \textup{7}RU, U.K.}
\email{Zinaida.Lykova@ncl.ac.uk}

\author{N. J.  Young}
\address{School of Mathematics, Statistics and Physics, Newcastle University, Newcastle upon Tyne NE1 7RU, U.K.
{\em and} School of Mathematics, Leeds University,  Leeds LS2 9JT, U.K.}
\email{Nicholas.Young@ncl.ac.uk}
\date{9th January 2019}

\begin{abstract} 
The symmetrized bidisc 
\[
G \stackrel{\rm{def}}{=}\{(z+w,zw):|z|<1,\ |w|<1\}
\]
has interesting geometric properties.  While it has a plentiful supply of complex geodesics and of automorphisms, there is nevertheless a unique complex geodesic $\calr$ in $G$ that is invariant under all automorphisms of $G$. 
Moreover, $G$ is foliated by those complex geodesics that meet $\calr$ in one point and have nontrivial stabilizer.
  We prove that these properties, together with two further geometric hypotheses on the action of the automorphism group of $G$, characterize the symmetrized bidisc in the class of complex manifolds.

\end{abstract} 

\maketitle
\tableofcontents
\section*{Introduction}

By a {\em domain} we mean a connected open set in $\c^n$ for some integer $n\geq 1$.  
\index{domain}
A domain is {\em homogeneous}
\index{homogeneous}
 if the automorphisms
 of the domain act transitively.  It is {\em symmetric}
\index{symmetric}
if every point of the domain is an isolated fixed point of an involutive automorphism of the domain.
\index{domain!bounded symmetric homogeneous}

The nature of a bounded symmetric homogeneous domain in $\c^n$ is captured by the great classification theorem of \'Elie Cartan \cite{ecartan},
\index{Cartan!\'E.}
an early triumph of the theory of several complex variables \cite{fuks,helg}.
\index{Cartan's theorem}
It states that any such domain is isomorphic to a product of domains, each of which is isomorphic to a domain of one of six concrete types.  The theorem is fundamental to the complex geometry and function theory of bounded symmetric homogeneous domains.

In this paper we are interested in  irreducible domains $\Omega$ which 
 narrowly miss being homogeneous, in the sense that the action of the automorphisms of $\Omega$ splits the domain into a one-parameter family of orbits.    
\index{orbit}
Such domains are said to have {\em cohomogeneity $1$}, and have an extensive theory \cite{GZ2,CGLP}  in both the mathematical and physics literatures.

One familiar domain that has cohomegeneity $1$ is the annulus
\index{annulus}
\index{$\mathbb{A}_q$}
\beq\label{defAq}
\mathbb{A}_q\df \{z\in\c: q<|z|<q\inv\}
\eeq
where $0<q<1$.  The orbits here are the sets
\beq\label{orbsA}
\{z: |z|=t\} \cup \{z:|z|=t\inv\}
\eeq
where $q<t\leq 1$.

For a higher-dimensional example,
 consider the domain
\beq\label{defG}
G\df\{(z+w,zw):|z|<1,\ |w|<1\}
\eeq
in $\C^2$, known as the {\em symmetrized bidisc}.  
\index{$G$}
\index{symmetrized bidisc}
The automorphisms of $G$ are the maps of the form
\beq\label{formautg}
(z+w,zw) \mapsto (m(z)+m(w),m(z)m(w))
\eeq
for some automorphism $m$ of the unit disc $\d$.  The orbits in $G$ are therefore generically $3$-dimensional real manifolds, and there is a one-parameter family of them.

Another domain, now in $\c^3$, having a one-parameter family of orbits is the {\em tetrablock}, 
\index{tetrablock}
which comprises the points $(x^1,x^2,x^3) \in\c^3$ such that 
\beq \label{deftet}
 1-x^1z-x^2w+x^3 zw\neq 0 
\eeq
for all  $z,w \in \c$   such that $|z| \le 1$ and $|w| \le1$. 

An ambitious project would be to classify bounded domains in $\c^n$, and more generally, complex manifolds, for which the orbits under the automorphisms of the manifold comprise a one-parameter family.  By way of a start we shall here characterize in geometric terms our archetypal example $G$ defined in equation \eqref{defG}.  This domain has been studied
 by numerous authors over the past 20 years, and has proved to be
a domain with a very rich complex geometry and function theory: see, besides many other papers, 
\cite{ay2004,cos04,ez05,jp04,pz05,kos,tryb,aly2016}.   $G$ is significant for the theory of invariant distances \cite{jp},
because it has Lempert's property, that 
the Carath\'eodory and Kobayashi metrics coincide
\cite{lem81}, despite the fact that $G$ is not convex (nor even biholomorphic to a convex domain  \cite{cos04}).
\index{Lempert}
It plays a role in operator theory \cite{bhatta,sarkar} and even has applications
to a problem in the theory of robust control (for example, \cite{Y12}); indeed the control application was the original motivation for the study of $G$.
In an earlier paper \cite{aly2018} we characterized $G$ in terms of the Carath\'eodory extremal functions that it admits.  Here we give another characterization,
this time in terms of its complex geodesics and  automorphisms.

 An {\em automorphism} of a complex manifold  is a bijective holomorphic self-map of the manifold; such a map automatically has a holomorphic inverse.
\index{automorphism}
For any complex manifold  $\Omega$ we denote by $\aut\Omega$ the automorphism group of $\Omega$ with the compact-open topology.
A {\em complex geodesic} of $G$ can be defined as the range of an analytic map $f:\D\to G$ that has an analytic left inverse, where $\D=\{z\in\c:|z|< 1\}$.
\index{$\D$}
\index{geodesic!complex}

We draw attention to two striking geometric properties of $G$.
\begin{enumerate}
\item  There exists a unique complex geodesic $\calr$ in $G$  that is invariant under all automorphisms of $G$.
Moreover,
every automorphism of $\calr$ extends to a unique automorphism of $G$;
\item  
for every $s\in\calr$ there exists a unique geodesic $F_s$ in $G$  having a nontrivial stabilizer in $\aut G$ and such that
\[
F_s\cap\calr = \{s\}.
\]
 Moreover, the geodesics $\{F_s:s\in\calr\}$ foliate $G$.\footnote{That is, every point of $G$ lies in some $F_s$ and no point of $G$ lies in two distinct $F_s$}
\end{enumerate}
We call $\calr$ the {\em royal variety} and the sets $F_s$ the {\em flat geodesics} of $G$.  
\index{royal!variety}
\index{geodesic!flat}

Could it be that properties (1) and (2) suffice to characterize $G$?  In the present paper we show that the answer is yes under some further geometric hypotheses,
 which we now describe.

We say that a properly embedded analytic disc\footnote{not assumed to be a geodesic} $D$ in a complex manifold $\Omega$
is a {\em royal disc} if it has properties analogous to those of $\calr$ in (1), 
\index{royal!disc}
that is,
$D$ is invariant under every automorphism of $\Omega$, and every automorphism of $D$ extends to a unique automorphism of $\Omega$.  
 A {\em royal manifold} is a pair $(\Omega, D)$  where $\Omega$ is a complex manifold and $D$ is a royal disc in $\Omega$.

If $(\Omega, D)$ is a royal manifold then a collection $\cale=\{E_\la:\la\in D\}$ of properly embedded analytic discs\footnote{again, not assumed to be geodesics} in $\Omega$ is a {\em flat fibration over $D$} if it has properties similar to those of $\{F_s:s\in \calr\}$ in (2), that is, $E_\la\cap D=\{\la\}$ for every $\la\in D$,  $\cale$ is a partition of $\Omega$ and, for every automorphism $\theta$ of $\Omega$ and every $\la\in D$,  $\theta(E_\la)= E_{\theta(\la)}$.  
 The triple $(\Omega,D, \cale)$ is then called a {\em flatly fibered royal manifold}.

The orbits in $(\Omega,D,\cale)$ have a natural parametrization by $[0,\infty)$.
For any $\mu\in\Omega$ there is a unique $\la\in D$ such that $\mu\in E_\la$;
we  define  the {\em Poincar\'e parameter} $P(\mu)$ 
\index{Poincar\'e parameter}
to be  the Poincar\'e distance from $\mu$ to $\la$ in the disc $E_\la$ (see Definition \ref{defP}).
Two points $\mu_1,\mu_2$ in $\Omega$ lie in the same orbit if and only if $P(\mu_1)=P(\mu_2)$.

Flatly fibered royal manifolds
 can enjoy two geometric properties: synchrony and sharpness. 
 {\em Synchrony} 
\index{synchrony}
is a condition which relates the actions of 
$\aut\Omega$ on $D$ and on the discs in $\cale$.  To be precise, if $\theta$ is an automorphism of $\Omega$ which fixes a point $\la\in D$, then it follows easily from the definition of a flat fibration over a royal manifold that the eigenspaces of the operator $\theta'(\la)$ on the tangent space $T_\la\Omega$ to $\Omega$ at $\la$
\index{$T_\la\Omega$}
are the tangent spaces $T_\la D$ and $T_\la E_\la$.  We say that $(\Omega,D,\cale)$ is {\em synchronous}  if, for every $\la\in D$,  
the eigenvalue of $\theta'(\la)$ corresponding to $T_\la E_\la$ is the square of the eigenvalue of $\theta'(\la)$ corresponding to $T_\la D$.
 
{\em Sharpness} 
\index{sharp!action}
is a condition on the action of 
$\aut\d$ on $\Omega$ in a flatly fibered royal manifold $(\Omega,D,\cale)$.
The definition of $(\Omega,D,\cale)$ implies that every  $m\in\aut\d$ induces an automorphism $\Theta(m)$ of $\Omega$.
For $\al\in\d$,  let $B_\al$ denote the automorphism
\beq\label{defBal}
B_\al(z)=\frac{z-\al}{1-\bar\al z}
\eeq
of $\d$.  
\index{$B_\al$}
We say that {\em $\aut\Omega$ acts sharply} at a point $\mu\in\Omega\setminus D$ if, in local co-ordinates,
\beq\label{sharp-act}
\e^{2P(\mu)}(\Theta(B_{it})(\mu) -\mu) =i\left(\Theta(B_t)(\mu)-\mu\right) + o(t)
\eeq
as $t \to 0$ in $\r$.
\black

The geometric content of the sharpness condition at $\mu$ relates to the derivative at zero of the map $\al\mapsto \Theta(B_\al)(\mu)$ from $\d$ to $\Omega$.  This map is {\em a priori} a real-linear map from $T_0\d$ to $T_\mu\Omega$; now $T_0\d(=\c)$ and $T_\mu\Omega$ are both complex vector spaces, and the sharpness condition is equivalent to the statement that the derivative at zero is also a {\em complex} linear map.

If we denote by $\mu^\sharp$ the range of this complex-linear derivative, then it is easy to see that $\mu^\sharp$ 
\index{$\la^\sharp$}
is the unique nonzero complex linear subspace of  $T_\mu\Omega$ that is
 contained in the $3$-dimensional real tangent space at $\mu$ to the orbit of $\mu$ in $\Omega$.  

The sharp direction
\index{sharp!direction}
 $\mu^\sharp$ is a covariant line bundle over $\Omega$ which has interesting geometric properties.  For example, in $G$, the sharp direction $\mu^\sharp$ is characterized by the fact that the complex geodesic $\calc$ through $\mu$ with direction $\mu^\sharp$ has the {\em closest point property},
meaning that, for any point $\la\in\calc$, if $F_s, \, s\in\calr$, is the flat geodesic containing 
$\mu$, then the closest point to $\la$ in $F_s$ is $\mu$.

Our main result, Theorem \ref{autprop70} in the body of the paper, gives a precise version of the following statement, which holds under suitable regularity conditions.\\
\vspace*{0.15cm}

{\bf Theorem A.} \em
Let $\Omega$ be a complex manifold.  $\Omega$ is isomorphic to $G$ if and only if there exist
a royal disc $D$ in $\Omega$ and a flat fibration $\cale$ of $\Omega$ over $D$ such that 
 $(\Omega,D,\cale)$ is a synchronous flatly fibered royal manifold 
 and $\aut{\Omega}$ acts sharply on $E_\la\setminus\{\la\}$ for some $\la\in D$.

\vspace*{0.5cm}
\rm

Formal definitions of synchrony and sharp action are given in Subsections \ref{synchr}
and \ref{sharpness}.  The appropriate notion of regularity is described in Subsection \ref{regroyal}.

Remarkably,  Theorem A implies that if $(\Omega,D,\cale)$ is a synchronous flatly fibered royal domain with suitable regularity, and $\aut\Omega$ acts sharply, then both $D$ and the leaves in $\cale$ are complex geodesics of $\Omega$. 
 It suggests that $G$ might be characterized also in terms of the properties of its complex geodesics, and in a future paper we shall show that it is so. 

In Section \ref{charactflat} we give  in Theorem \ref{thm30} a characterization of $G$ 
 in terms of the existence of global co-ordinates ranging over the bidisc and satisfying certain partial differential equations.  These co-ordinates are related to the flat geodesics in $G$.

\black
In a short final section we discuss the relevance of the notion of symmetric space to the question of classification and show that the annulus, the symmetrized bidisc and the tetrablock, besides being inhomogeneous, also fail to be symmetric in \'E. Cartan's sense.

If $U$ and $\Omega$ are complex manifolds, we denote by $\Omega(U)$ the set of holomorphic mappings from $U$ into $\Omega$. 
\index{$\Omega(U)$}

We have used the expression {\em properly embedded analytic disc} in a complex manifold $\Omega$.  By this phrase we mean a proper injective analytic map $k:\d\to\Omega$ such that 
$k'(z)\neq 0$ for all $z\in \d$.  The range of such a map $k$ will also be called a properly embedded analytic disc.
\index{analytic disc!properly embedded}

\section{The action of automorphisms on $G$} \label{actionG}
 In this section we study the orbit structure of $G$ under the action of $\aut G$.

\subsection{The action of $\aut \d$ on $G$}

As we stated in the introduction (see equation \eqref{formautg}), every automorphism $m$ of $\d$ induces a map $\gamma_m:G \to G$ via the formula
\be\label{aut10}
\gamma_m (z+w,zw) = (m(z)+m(w),m(z)m(w))
\ee
for $z,w \in\d$.
\index{$\ga_m$}
It is easy to check that this formula defines a map $\gamma_m \in G(G)$ and that $\gamma_m \in \aut  G$.
\begin{prop} \label{gamiso}
 The map $\ga:\aut \d \to \aut  G$ given by
\be\label{aut20}
 \ga(m) = \gamma_m 
\ee
for $m\in\aut \d$ is a continuous isomorphism of topological groups.
\end{prop}

 The fact that $\ga$ is an isomorphism of groups is proved in \cite[Theorem 5.1]{AY08} or \cite{jp}.  
It is routine to show that $\ga$ is continuous with respect to the compact-open topologies on $\aut\d$ and $\aut G$.

The following statements are elementary.
\begin{prop}\label{diffeval}
\begin{enumerate}[\rm (1)]
\item $\aut\d$ and $\aut G$ are Lie groups.
\item For any $s\in G$ the map 
\beq\label{defes}
e_s:\aut\d\to G \mbox{  given by } e_s(m)=\ga_m(s)
\eeq
 is real-analytic.
\end{enumerate}
\end{prop}
\index{Lie!group}
The map $e_s$, where $s\in G$, will be called the {\em evaluation map} at $s$ on $\aut\d$.
\index{evaluation map}
\index{$e_s$}

\subsection{The action of $\aut G$ on the royal variety}
The {\em royal variety} in $G$ is defined to be the set
\begin{align*}
\calr  &= \{s=(s^1,s^2)\in G: (s^1)^2=4s^2\} \\
	&=\{(2z,z^2):z\in\d\}
\end{align*}
(we use superscripts to denote the components of a point in $\C^d$).
Thus  $\calr=R(\d)$ where 
\be\label{aut42}
R(z) = (2z,z^2)\qquad \mbox{ for } z\in \d.
\ee
Clearly
\be\label{aut44}
\gamma_m (R(z)) = R( m (z))\qquad \mbox{ for } z\in \d  \mbox{ and }m\in \aut \d.
\ee

The observations \eqref{aut20} and \eqref{aut44} have three consequences, summarized in the following proposition.
\begin{prop}\label{autprop10}
\begin{enumerate}[\rm (1)]
\item Every automorphism of $G$ leaves $\calr$ invariant.
\item Every automorphism of $G$ is uniquely determined by its values on $\calr$.
\item Every automorphism of $\calr$ has a unique extension to an automorphism of $G$.
\end{enumerate}
\end{prop}

In statement (3), automorphisms of $\calr$ are with respect to the structure of $\calr$ as a complex manifold.

We can summarize these three statements by saying that the restriction map $\ga \mapsto \ga|\calr$ is an isomorphism from $\aut  G$ to $\aut\calr$.  The following commutative diagram describes the situation, where $\iota_\calr$ denotes the injection of $\calr$ into $G$ and $m\in\aut \d$.
\be\label{commDiagG}
\begin{array}{rrrcl}
\d & \stackrel{R}{\longrightarrow} &\calr& \stackrel{\iota_\calr}{\longrightarrow} & G \\
\vcenter{\llap{$\scriptstyle{m}$}}\Big\downarrow&  &\vcenter{\llap{$\scriptstyle{\ga_m|\calr}$}}\Big\downarrow & & \vcenter{\llap{$\scriptstyle{\ga_m}$}}\Big\downarrow\\
\d & \stackrel{R}{\longrightarrow} &\calr& \stackrel{\iota_\calr}{\longrightarrow} & G
\end{array}
\ee

\subsection{Orbits in $G$ as manifolds}
For any complex manifold $U$ and any $\lambda \in U$,  we denote by $\orb_U(\lambda)$
\index{$\orb_U(\lambda)$}
the orbit of $\lambda$ under the action of the group of automorphisms of $U$:
\[
\orb_U(\lambda) = \{\gamma(\lambda): \gamma \in \aut U\}.
\]

Consider the case that $U=G$ and $\lambda=s\in G$.  

In view of Proposition \ref{gamiso}, for any $s\in G$,
\beq\label{formOrbs}
\orb_G(s)=\{\ga_m(s): m\in\aut\d\},
\eeq
so that  $\orb_G(s)$ is the range of the evaluation map  $e_s$ of equation \eqref{defes}.

$\aut\d$ is a $3$-dimensional  real-analytic manifold, for which we shall need local co-ordinates.
\begin{lemma}\label{atlas}
For $(r,\alpha) \in \r \times \d$ let $m_{r,\alpha} \in \aut\d$ be given by the formula
\be\label{aut90}
m_{r,\alpha} (z) =\e^{ir}\frac{z-\alpha}{1-\bar\alpha z}, \qquad z\in\d.
\ee
Let
\begin{align}\label{U1U2}
U_1&=\{m_{r,\alpha}: -\pi < r<\pi, \, \al\in\d\}, \\
   U_2&=\{m_{r,\alpha}: 0 < r<2\pi, \, \al\in\d\} \notag
\end{align}
and define
\[
\ph_1:U_1 \to (-\pi,\pi)\times \d \quad \mbox{ by } \quad  \ph_1(m_{r,\alpha})=(r, \al),
\]
and similarly for $\ph_2:U_2\to (0,2\pi)\times \d$.
 Then $(U_1,\ph_1)$ and $(U_2,\ph_2)$ are charts in $\aut\d$ which together comprise a real-analytic atlas $\cala$ for the group manifold.  The identity automorphism $\idd=m_{0,0}$ belongs to $U_1$.  
\end{lemma}

\begin{proof}
The automorphisms of $\d$ consist of the maps $m_{r,\al}$ for $r\in [-2\pi,2\pi]$ and $\al\in\d$, and therefore 
 $\aut\d=U_1 \cup U_2$.  If $-\pi < r < 0$ then
\[
\ph_2\circ \ph_1\inv (r,\al)= (r+2\pi,\al)
\]
and similarly when $0<r<\pi$. 
The transition map is therefore real-analytic from $\ph_1(U_1\cap U_2)$ to $\ph_2(U_1\cap U_2)$.
 \end{proof}

\begin{prop}\label{autrem15} 
If $s \not\in \calr$ then the evaluation map  $e_s:\aut\d\to \orb_G(s)$
 is a local homeomorphism and a two-to-one covering map, given explicitly by
\[
e_s(m_{r,\al})= \frac{\left(\e^{ir}(-2\al+(1+|\al|^2)s^1-2\bar\al s^2), \e^{2ir}(\al\al-\al s^1+s^2)\right)}{1-\bar\al s^1+\bar\al\bar\al s^2}.
\]
\end{prop}

\begin{proof}
Consider a point $s=(z+w,zw)\in G$ where $z,w \in \d$ and $z \neq w$.   Let $\up$ be the unique automorphism of $\d$ that maps 
$z$ to $w$ and $w$ to $z$.  Note that  $ \up$ 
is not the identity automorphism ${\mathrm{id}}_{\d}$ since $z \neq w$.
For $m_1,m_2 \in \aut\d$,
\begin{align}
\ga_{m_1}(s)=\ga_{m_2}(s) & \Leftrightarrow \pi(m_1(z),m_1(w))=\pi(m_2(z),m_2(w)) \notag \\
	&  \Leftrightarrow  \left\{ \begin{array}{l}\mbox{ either } m_1(z)=m_2(z) \mbox{ and } m_1(w)=m_2(w) \\
						\mbox{ or } m_1(z)=m_2(w) \mbox{ and } m_1(w)=m_2(z)  \end{array} \right. \notag\\
	&  \Leftrightarrow  \left\{ \begin{array}{l}\mbox{ either } m_1=m_2 \\
						\mbox{ or } m_2\inv\circ m_1=\up.
						   \end{array} \right. \label{2to1}
\end{align}
  Thus $e_s:m \mapsto \ga_m(s)$ is two-to-one from $\aut\d$ to $ \orb_G (s)$.

To prove that $e_s$ is a local homeomorphism, choose any point $e_s(\beta)$ of $ \orb_G (s)$, where $\beta\in \aut\d$.
Choose a neighborhood $U$ of $\idd$ such that 
\beq\label{m1m2}
m_2\inv \circ m_1\neq \up\quad \mbox{  for all } m_1, m_2 \in U.
\eeq
We claim that $\beta\circ U$
 is a neighborhood of $\beta$ on which $e_s$ 
is injective.  Certainly it is a neighborhood of $\beta$, and if $e_s|\beta \circ U$ is not injective then there exist distinct points $m_1,m_2 \in U$ such that $e_s(\beta\circ {m_1})=e_s(\beta\circ {m_2})$.  That is, $\ga_\beta\circ\ga_{m_1}(s)=\ga_\beta\circ \ga_{m_2}(s)$, and therefore $\ga_{m_1}(s)=\ga_{m_2}(s)$.
Hence, by the equivalence \eqref{2to1}, $m_2\inv \circ m_1=\up$.  This equation contradicts the statement \eqref{m1m2}.  Thus $e_s$ is locally injective on $\aut{\d}$.

Choose a compact neighborhood $V$ of $\idd$ contained in $U$.  Since a continuous bijective map from a compact space to a Hausdorff space is a homeomorphism, $e_s|V$ is a homeomorphism onto its range.  It follows by homogeneity that $e_s$ is a local homeomorphism on $\aut\d$.  Indeed, consider any $m\in \aut\d$ and its neighborhood $m\circ V$.   Define $L_m:\aut\d\to\aut\d$ by $L_m(\theta)= m\inv\circ\theta$ for $\theta\in\aut\d$. In the commutative diagram
\be\label{commDiagFs}
\begin{array}{rrrcl}
m\circ V & \stackrel{e_s|m\circ V}{\longrightarrow} &e_s(m\circ V)\\
\vcenter{\llap{$\scriptstyle{L_m|m\circ V}$}}\Big\downarrow&  &\vcenter{\llap{$\scriptstyle{\ga(m)|e_s(V)}$}}\Big\uparrow\\
V & \stackrel{e_s|V}{\longrightarrow} & e_s(V) 
\end{array}
\ee
the map $e_s|m\circ V$ is expressed as the composition of three homeomorphisms, and so is itself a homeomorphism.  Thus $e_s$ is a local homeomorphism.

The formula for $e_s(m_{r,\al})$ is a simple calculation.
\end{proof}
 For any $s\in G$, the map $e_s'(\idd)$ is a real-linear map from the tangent space $T_{\idd}\aut\d$ to $T_s\orb_G(s)$.  The space $T_{\idd}\aut\d$ is the Lie algebra of $\aut\d$, so we shall denote it by $\laut$
\index{$\laut$}
\index{Lie!algebra}
 (though we shall not use its Lie structure, only its real-linear structure).  

For every $s\in G$ we define a real-linear subspace $\calv(s)$
\index{$\calv(s)$}
 of $\C^2$ by
\be\label{basisV}
\calv(s) \df \spa_\r\left\{ i\bpm s^1\\2s^2\epm, \bpm 2-(s^1)^2 + 2s^2 \\s^1-s^1s^2 \epm, i\bpm 2 +(s^1)^2 -2s^2 \\ s^1 + s^1s^2 \epm  \right\}.
\ee
\begin{thm}\label{autprop30a}
\begin{enumerate}[\rm (1)]
\item If $s \in \calr$, then $\orb_G (s)$ is a one-dimensional complex manifold properly embedded in $G$.  

\item If $s \in G \setminus  \calr$, then 
$\orb_G (s)$ is a three-dimensional real-analytic manifold properly embedded in $G$.   
\end{enumerate}
Moreover, in either case, the tangent space to $\orb_G(s)$ at $s$ is $\calv(s)$ and
\beq\label{calvran}
\calv(s)=\ran e_s'(\idd).
\eeq
\end{thm}
In the sequel the notation $T_s\orb_G(s)$ denotes the complex tangent space if $s\in\calr$ and the real tangent space if $s\notin\calr$.
Thus, for all $s\in G$,
\beq\label{Tsranes}
T_s\orb_G(s)= \ran e_s'(\idd).
\eeq
\begin{proof}

Consider $s\in G$. 
We shall calculate the rank of the real linear operator $e_s'(m)$ for $m\in\aut \d$.  Let $\iota : \orb_G(s) \to \C^2_r$ denote the inclusion map.

\begin{lemma}\label{vralpha}
For any tangent vector $(r,\al)$ at $(0,0)$ to $(-\pi,\pi)\times\d$, let $v_{r,\al}(s)$ denote the tangent vector $(e_s\circ\ph_1\inv)'(0,0)(r,\al)$
 in $T_sG \subset \C^2$.  Then
\beq \label{formvral}
v_{r,\al}(s)=ir\begin{pmatrix}s^1\\2s^2\end{pmatrix}
	-\alpha\begin{pmatrix}2\\s^1\end{pmatrix}
	+\bar\alpha\begin{pmatrix}{(s^1)}^2-2s^2\\s^1s^2\end{pmatrix}.
\eeq
\index{$v_{r,\al}$}
\end{lemma}
\begin{proof}
We have $(r,\al)\in \r\times\c$.  Define a path $\kappa(t)=(tr,t\al)$ in $(-\pi,\pi)\times\d$ for $|t|<\eps$, where $\eps$ is small enough.

Then let 
\[
f_s=\iota \circ e_s\circ \ph_1\inv :(-1,1)\times \D \to \C^2_r
\]
 and define $v_{r,\alpha} (s) \in \c^2$ by the formula
\begin{align}
v_{r,\alpha}(s)&=\frac{\dd}{\dd t} f_s \circ\kappa(t)\big|_{t=0}\label{defvra}  \\
	&=\frac{\dd}{\dd t} \iota\circ e_s\circ \ph_1\inv\circ\kappa(t)\big|_{t=0}\notag  \\
	&= \frac{\dd}{\dd t} \iota\circ e_s\circ m_{\kappa(t)}\big|_{t=0} \notag \\
	&=\frac{\dd}{\dd t}\iota\circ e_s\circ m_{tr,t\al}\big|_{t=0} \notag \\
	&= \frac{\dd}{\dd t}\iota\circ  \ga_{ m_{tr,t\al}}(s)\big|_{t=0}.\label{aut92}
\end{align}
From equation \eqref{aut90}, for any $z\in\d$, 
\[
\frac{\dd}{\dd t} m_{tr,t\alpha} (z) \big|_{t=0} = irz - \alpha +\bar \alpha z^2.
\]

Hence, by equations \eqref{aut10} and \eqref{aut92}, if $s=(z_1+z_2,z_1z_2)$,
\begin {align}\label{aut95}
v_{r,\alpha} (s) &=  \frac{\dd}{\dd t} \left. \bpm m_{tr,t\al}(z_1)+m_{tr,t\al}(z_2) \\ m_{tr,t\al}(z_1) m_{tr,t\al}(z_2) \epm\right|_{t=0} \notag \\
	&= \bpm irz_1-\al+\bar\al(z_1)^2 + irz_2-\al+\bar\al (z_2)^2 \\
		(  irz_1-\al+\bar\al(z_1)^2)z_2 + (irz_2-\al+\bar\al (z_2)^2)z_1 \epm  \notag \\
	&=ir\begin{pmatrix}s^1\\2s^2\end{pmatrix}
	-\alpha\begin{pmatrix}2\\s^1\end{pmatrix}
	+\bar\alpha\begin{pmatrix}{(s^1)}^2-2s^2\\s^1s^2\end{pmatrix}.
\end{align}
\end{proof}
Continuing the proof of Theorem \ref{autprop30a}, by the Chain Rule we have, from equation \eqref{defvra},
\begin{align}
v_{r,\al}(s) &= f_s'(\kappa(0))\kappa'(0) \notag \\
	&= f_s'(0,0) \bpm r\\  \al\epm. \label{formvra2}
\end{align}
Thus the range of the real linear map $f_s'(0,0): \r\times\c \to \c^2$ is the set
\begin{align*}
\ran f_s'(0,0) &= \{v_{r,\al}: r\in\r, \al\in\c\} \\
	&=\left\{ir\begin{pmatrix}s^1\\2s^2\end{pmatrix}
	-\alpha\begin{pmatrix}2\\s^1\end{pmatrix}
	+\bar\alpha\begin{pmatrix}{(s^1)}^2-2s^2\\s^1s^2\end{pmatrix}:r\in\r, \, \al\in\c\right\}.
\end{align*}
On taking $(r,\al)$ to be successively $(1,0), \, (0,-1)$ and $(0,-i)$ we find that, for any $s\in G$,
\[
\ran f_s'(0,0) = \calv(s),
\]
the real vector space introduced in equation \eqref{basisV}.  Thus
\be\label{rankFs}
\iota'(s)\ran e_s'(\idd)=\ran f_s'(0,0) = \calv(s)
\ee
for all $s\in G$.    In the sequel we shall suppress the inclusion map $\iota'(s)$ and regard $\ran e_s'(\idd)$ as a subspace of $\C^2_r$.

Now consider $s\in G\setminus \calr$.  
By Lemma \ref{autlem8} below, $\dim_\r \calv(s) = 3$, and so, by equation \eqref{rankFs},
$e_s'(\idd)$ has rank $3$.  We claim that $e_s'(m)$ has rank $3$ for all $m\in \aut \d$.   Indeed, on differentiating the relation
\[
e_s(m)=\ga_m(s)=\ga_m\circ \ga_{\idd}(s)=\ga_m\circ e_s(\idd),
\]
we find (since $e_s(\idd)=s$) that
\[
e_s'(m)= \ga_m'(s)e_s'(\idd).
\]
Since $\ga_m$ is an automorphism of $G$, $\ga_m'(s)$ is a nonsingular real linear transformation of $\c^2$.   Thus
\[
\rank_\r e_s'(m) = \rank_\r e_s'(\idd) =3
\]
for every $m\in\aut\d$.

We wish to deduce that $\orb_G(s)$ is a real $3$-dimensional $C^\infty$-manifold which (modulo the identification map $\iota'(s)$) lies in $\c^2$.
The following statement is \cite[Theorem 5.2]{Spiv}.

 {\em A subset $M$ of $\r^n$ is a $k$-dimensional manifold if and only if, for every point $s\in M$ there exist an open neighborhood $V$ of $s$ in $\r^n$, an open set $W$ in $\r^k$ and an injective differentiable function $f:W \to \r^n$ such that 
\begin{enumerate}[\rm (1)]
\item $f(W)=M\cap V$,
\item $f'(y)$ has rank $k$ for every $y\in W$.
\end{enumerate}
}

We shall apply this criterion in the case $n=4, \ k=3, \ M=\orb_G(s)$.    Consider any point $e_s(m)\in\orb_G(s)$, where $m\in\aut\d$, say $m\in U_j, \, j=1$ or $2$.  By Proposition \ref{autrem15}, $e_s$ is a local homeomorphism, and so we may choose an open neighborhood $N$ of $m$ in $U_j$ such that $e_s|N$ is a homeomorphism from $N$ to an open subset of $\orb_G(s)$.  Since $\orb_G(s)$ has the relative topology induced by $G$, there is an open set $V$ in $\r^4$ such that $e_s(N)= V\cap \orb_G(s)$.

Let $W=\ph_j(N)$.  Then the map $f= e_s\circ \ph_j\inv$ satisfies conditions (1) and (2).  It follows that $\orb_G(s)$ is a real $3$-dimensional $C^\infty$ manifold in $\c^2$.

The linear map $(e_s\circ \ph_1\inv)'(0,0)$ maps  the tangent space $T_{(0,0)} (-\pi,\pi)\times\d$ into the tangent space $T_s\orb_G(s)\subset \c^2$.    We have seen that the range of $(e_s\circ \ph_1\inv)'(0,0)$ is $\calv(s)$.  Hence 
\[
\calv(s) \subseteq  T_s\orb_G(s).
\]
Since both spaces have real dimension $3$, the inclusion holds with equality.

In the case that $s\in\calr$, say $s=(2\zeta,\zeta^2)$ for some $\zeta\in\d$,
\[
\orb_G(s)=\calr=\{(2z,z^2): z\in\d\},
\]
 which is  a one-dimensional complex manifold properly embedded in $G$ by the map $R:\d\to G$.
The complex tangent space to $\calr$ at $s$ is $\c(1,\zeta)$, and, by equation \eqref{basisV},
\begin{align*}
\calv(s)&=\spa_\r\left\{ i\zeta\bpm 1 \\ \zeta\epm, (1-\zeta^2)\bpm 1 \\ \zeta \epm, i(1+\zeta^2)\bpm 1 \\ \zeta \epm \right\} \\
	&=\c\bpm 1 \\ \zeta \epm.
\end{align*}
Thus $T_s\orb_G(s) = \calv(s)$ in the sense of complex manifolds.
\end{proof}

\begin{lem}\label{autlem8}
For any $s\in G$, the real vector space $\calv(s)$ defined by equation \eqref{basisV} satisfies
\be\label{aut100}
\dim_\r  \calv(s) =\left\{ \begin{array}{lll}3 &\quad &\mbox{ if }  (s^1)^2\neq 4 s^2  \\  
								2    &\quad & \mbox{ if }  (s^1)^2 = 4 s^2.  \end{array} 
			    \right.
 \ee
\end{lem}

\begin{proof}

It is clear from the definition \eqref{basisV} that $\calv(s)$ is a real vector subspace of $\c^2$ of real dimension at most $3$.

Suppose that scalars $\la,\mu,\nu \in\r$ satisfy
\be\label{linrel}
\la i\bpm s^1\\2s^2\epm + \mu \bpm 2-(s^1)^2 + 2s^2 \\s^1-s^1s^2 \epm + \nu i\bpm 2 +(s^1)^2 -2s^2 \\ s^1 + s^1s^2 \epm =0.
\ee
Multiply on the left by the row matrix $\bpm 2s^2 & -s^1 \epm$ to obtain
\be\label{elimla1}
-((s^1)^2-4s^2)\left((1+s^2)\mu + (1-s^2)\nu i \right)=0.
\ee
Consider the first case in equation \eqref{aut100}, namely, that  $(s^1)^2\neq 4s^2$ (equivalently,  $s\notin \calr$).  By equation \eqref{elimla1}
\[
(1+s^2)\mu + (1-s^2)\nu i=0,
\]
whence
\[
\mu+i\nu= -s^2(\mu-i\nu).
\]
Since $s\in G$, we have $|s^2|<1$, and so necessarily $\mu=\nu=0$.  Since $(s^1)^2\neq 4s^2$, at least one of
$s^1,s^2$ is nonzero, and so, by equation \eqref{linrel}, $\la=0$.  Hence the three spanning vectors for $\calv(s)$ in equation \eqref{basisV} are linearly independent.  We have shown that $\dim_\r \calv(s) =3$ when $s\notin \calr$.

Next consider a point $s\in \calr$.   On substituting $s^2= \tfrac 14 (s^1)^2$ in equation \eqref{basisV} we obtain
\[
\calv(s)=\spa_\r\left\{ is^1\bpm 2\\ s^1 \epm, (4-(s^1)^2)\bpm 2\\ s^1 \epm, i(4+(s^1)^2) \bpm 2\\ s^1 \epm    \right\}.
\]
Since each of these vectors is a complex scalar multiple of the vector $\bpm 2 & s^1 \epm^T$, it follows that $\dim_\r \calv(s) \leq 2$.

In fact $\dim_\r \calv(s) = 2$.  For otherwise the second and third spanning vectors for $\calv(s)$ are linearly dependent over $\r$, and so there exist $\mu,\nu \in\r$, not both zero, such that
\[
\mu ( 4-(s^1)^2)+\nu i (4+(s^1)^2) = 0
\]
and consequently
\[
4(\mu +\nu i)=(s^1)^2(\mu-\nu i).
\]
Thus $|s^1|=2$, contrary to choice of $s\in G$.
Therefore
\[
\dim_\r \calv(s) = 2 \quad \mbox{ when } s \in\calr.
\]
\end{proof}

\subsection{The sharp direction in $G$} \label{sharpinG}
By Theorem \ref{autprop30a}, for any  $s\in G\setminus \calr$, the tangent space  $\calv(s)$ at $s$ to the orbit $\orb_G(s)$ is a real $3$-dimensional subspace of $\c^2$.  Accordingly $\calv(s)$ contains a unique $2$-real-dimensional subspace that is also a one-dimensional complex subspace of $\C^2$, equal to $\calv(s)\cap i\calv(s)$.  On the other hand, for $s\in\calr$, the tangent space $\calv(s)$ is already a complex subspace of $\C^2$.
\begin{defin}\label{defsharp}
For any $s\in G$, the {\em sharp direction at $s$}
\index{sharp!direction}
\index{$s^\sharp$}
 is the unique nonzero complex subspace of $\calv(s)$ in $\C^2$ and is denoted by $s^\sharp$.  Thus
\[
s^\sharp = \calv(s) \cap i\calv(s).
\]
\end{defin}

The sharp direction is covariant with automorphisms of $G$, in the following sense.
\begin{prop}\label{autprop40}
If $\ga\in \aut  G$ and $s\in G$, then 
\[
\ga (s)^\sharp=\ga'(s) s^\sharp.
\]
\end{prop}
\begin{proof}
Since $\ga$ is a differentiable self-map of $\orb_G(s)$, its derivative $\ga'(s)$ is a real-linear map between the tangent spaces
$\calv(s)$ and $\calv(\ga(s))$.  Since furthermore $\ga'(s)$ is a nonsingular complex linear map from $T_s G\sim\C^2$ to $T_{\ga(s)}\sim \C^2$,
it maps  the complex subspace $s^\sharp$ of $\C^2$ to a nonzero complex subspace of $\C^2$.  Hence $\ga'(s)s^\sharp$ is a nonzero complex subspace of $\calv(\ga(s))$.  Hence $\ga'(s)s^\sharp=\ga(s)^\sharp$.
\end{proof}

\subsection{Flat geodesics and the action of $\aut  G$}\label{flatfibG}

In the introduction we defined the flat geodesics of $G$ to be the geodesics that meet the royal geodesic $\calr$ exacly once and are stabilized by a nontrivial automorphism of $G$.
\index{flat!geodesic}
This definition has the merit that it is geometrical in character, but in practice
(for example, to show that the flat geodesics foliate $G$) it is often simpler to use
the fact that  the flat geodesics in $G$ are the sets of the form
\beq\label{flatbeta}
F^\beta \df \{(\beta+\bar\beta z,z): z\in\d\}
\eeq
for some $\beta\in\d$. 
\index{$F^\beta$}
One can check that the point $s\in G$ lies on the unique $F^\beta$ with
\beq\label{forbeta}
\beta=\frac{s^1-\overline{ s^1}s^2}{1-|s^2|^2} \in \d.
\eeq
More details can be found in \cite{jp,cos04,AY06} and \cite[Appendix A]{aly2016}.

Let us at least sketch a proof that the set $F^\beta$ is indeed a flat geodesic according to the definition in the introduction.
Firstly, a straightforward calculation shows that any automorphism of $G$ maps $F^\beta$ to a set of the form $F^{\beta'}$ for some $\beta'\in\d$.
Clearly $F^\beta$ is a complex geodesic in $G$:  for any $\beta\in\d$  the co-ordinate function $s^2$ is a holomorphic left inverse of the properly embedded analytic disc $z\mapsto (\beta+\bar\beta z,z)$ in $G$. 
It is simple to check that $F^\beta$ meets $\calr$ exactly once, say at the point
$s(\beta) \in \calr$.  Choose a nontrivial automorphism $\theta$ of the analytic disc $\calr$ that fixes $s(\beta)$, and let $\ga$ be the unique extension of $\theta$ to a (necessarily nontrivial) automorphism of $G$.  Then $\ga(F^\beta)= F^{\beta'}$, and since $\ga$ fixes $s(\beta)$, it follows that $F^{\beta'} $ meets $\calr$ at $s(\beta)$.  Distinct sets $F^\beta$ are disjoint, and therefore $\beta=\beta'$.  That is, $F^\beta$ is stabilized by a nontrivial automorphism of $G$.

The converse statement, that every flat geodesic is an $F^\beta$, follows from the classification into five types of the complex geodesics in $G$ given in \cite[Chapter 7]{aly2016}.
\black

We summarize the main geometric properties of flat geodesics.
\begin{proposition}\label{mainflat}
\begin{enumerate}
\item[\rm{(1)}] Through each point $s$ in $G$ there passes a unique flat geodesic $F_s$.
\index{$F_s$}
\item[\rm{(2)}]  Every flat geodesic intersects the royal geodesic $\calr$ in exactly one point.
\item[\rm{(3)}]  Automorphisms of $G$ carry flat geodesics to flat geodesics.
\end{enumerate}
\end{proposition}
\noindent The following lemma is a reformulation of the first two of these facts.
\begin{lem}\label{autlem10}
The family
\[
\calf \df \{F_s:s \in \calr\}
\]
 is a partition of $G$.
\end{lem}
\begin{definition}\label{flatdir}
For any $s\in G$, the complex tangent space at $s$ to the unique flat geodesic through $s$ will be called the 
{\em flat direction} at $s$, and will be denoted by $s^\flat$.
\index{flat!direction}
\index{$s^\flat$}
\end{definition}
Thus, if $s^1= \beta+\bar\beta s^2$, then
\beq\label{formsflat}
s^\flat= \c(\bar\beta,1),
\eeq
which is a one-dimensional complex subspace of $\c^2$. 
 The map $s\mapsto s^\flat$ is a covariant line bundle which is a sub-bundle of $TG$.

Facts (1)-(3) in Proposition \ref{mainflat} imply the following description of the action of $\aut  G$ on $\calf$.
\begin{lem}\label{autlem20}
If $\gamma \in \aut  G$ and $s\in \calr$, then $\gamma (F_s) = F_{\gamma(s)}$.
\end{lem}
\begin{proof}
Fix $\gamma \in \aut  G$ and $s\in \calr$. By Fact 3, there exists $t\in \calr$ such that $\gamma(F_s) =F_t$, and Condition (i) in Proposition \ref{autprop10} implies that $\gamma(s) \in \calr$. Therefore 
$ \gamma(s) \in \calr\cap F_t$. Hence by Fact 2,  $t=\gamma(s)$.
\end{proof}
We shall call $\{F_s:s\in\calr\}$ the {\em flat fibration of $G$}.
\index{flat!fibration}

\begin{proposition}\label{flneqsh}
For all $s\in G$ the spaces $s^\sharp$ and $s^\flat$ are unequal.
\end{proposition}
This statement will follow from explicit formulae for the sharp and flat directions.
We already know that, for $s\in F^\beta$, $s^\flat$ is given by equation \eqref{formsflat}.
\begin{proposition}\label{formsharp}
For any $\beta\in\d$ and any $s\in F^\beta$,
\beq\label{formulasharp}
s^\sharp = \c \left(1, \frac {\beta-\half s^1}{1-\half\bar\beta s^1}\right).
\eeq
\end{proposition}
\begin{proof}
Consider $s\in G\setminus\calr$.  Let
\[
v_1= i\bpm s^1\\2s^2\epm, \quad v_2=\bpm 2-(s^1)^2 + 2s^2 \\s^1-s^1s^2 \epm, \quad v_3=i\bpm 2 +(s^1)^2 -2s^2 \\ s^1 + s^1s^2 \epm.
\]
By Theorem \ref{autprop30a}, $\{v_1,v_2,v_3\}$ constitutes a basis for $T_s\orb_G(s)$.
Let 
\[
c_1=-2s^1,\quad c_2=-i(1-s^2),\quad c_3=1+s^2
\]
and note that $c_2, c_3$ are nonzero.
One finds that $c_1v_1+c_2v_2+c_3v_3 =0$, and therefore
\[
(\re c_1)v_1+(\re c_2)v_2+(\re c_3) v_3=-i(\im c_1)v_1-i(\im c_2)v_2 -i(\im c_3) v_3.
\]
Hence
\[
v \df (\re c_1)v_1+(\re c_2)v_2+(\re c_3) v_3\neq 0
\]
and
\[
v \in T_s\orb_G(s) \cap iT_s\orb_G(s) = s^\sharp.
\]
Further calculation yields the formula
\[
v= 2i(1-|s^2|^2)\bpm 1-\bar\beta \half s^1 \\ \beta-\half s^1 \epm,
\]
in agreement with equation \eqref{formulasharp}.  This proves the proposition in the case that
$s\in G\setminus \calr$.  For $s=(2z,z^2)\in\calr$ equation \eqref{formulasharp} is easily checked.
\end{proof}
The fact that $s^\sharp \neq s^\flat$ can now be verified by a simple comparison of equations \eqref{formulasharp} and \eqref{formsflat}.

\begin{corollary}
The tangent bundle of $G$ is the direct sum of the sharp bundle and the flat bundle:
\[
T_sG= s^\sharp \oplus s^\flat \qquad \mbox{ for all }s\in G.
\]
\end{corollary}

\subsection{Synchrony in $G$}  \label{selfsynch}
\index{synchrony}
There is a subtle relationship between the action of an automorphism of $G$ on the royal variety and its action on any flat geodesic.

For any complex manifold $U$ and $\la$ in $U$,  denote by $ \aut_\la U$
\index{$\aut_\la$}
the stabilizer
\index{stabilizer}
 of $\la$ in $\aut U$ (also known as the isotropic subgroup of $\aut U$ at $\la$).
For any $s_0 \in \calr$, the sets $\calr$ and $F_{s_0}$ are embedded analytic discs in $G$ that intersect transversally at the point $s_0$. 
Every $\theta$ in $\aut_{s_0} G$ determines an automorphism of the analytic variety $\calr \cup F_{s_0}$. For an automorphism of a general variety there need be no connection between the action on  two leaves beyond what is implied by the condition that the restrictions of the automorphism to the two leaves must agree at any common point. However, in the context of the domain $G$, in the light of Condition (ii) in Proposition \ref{autprop10}, the action of $\theta$ on $F_{s_0}$ is uniquely determined by the action of $\theta$ on $\calr$. The following propositions describe this dependence explicitly.

We denote the unit circle $\{z:|z|=1\}$ in the complex plane by $\T$.
\index{$\T$}
For $\eta \in \t$ let $\rho_\eta$ denote the element of $\aut_0 \d$ defined by $\rho_\eta(z) = \eta z$. 
\index{$\rho_\eta$}
Clearly $\aut_0\d=\{\rho_\eta:\eta\in\t\}$.
\begin{proposition}\label{autlem30}
If $s_0\in \calr$ and $\theta\in \aut_{s_0} G$, then $\theta'(s_0)$ has eigenspaces $T_{s_0} \calr$ and $T_{s_0}F_{s_0}$ with corresponding eigenvalues $\eta$ and $\eta^2$ for some $\eta\in \t$.
\end{proposition}
\begin{proof}
Since $\theta$ leaves invariant both $\calr$ and $F_{s_0}$, it follows that $\theta'(s_0)$ leaves invariant
the tangent spaces $T_{s_0}\calr$ and $T_{s_0}F_{s_0}$.  These two one-dimensional tangent spaces are thus eigenspaces of $\theta'(s_0)$.

Observe that $\ga_{\rho_\eta}$ is the restriction to $G$ of the linear operator on $\c^2$ with matrix $\diag\{\eta,\eta^2\}$, and hence
\[
\ga_{\rho_\eta}'(s_0) \sim \diag\{\eta,\eta^2\}.
\]

Let $s_0=(2\al,\al^2)$. Since $\theta\in \aut_{s_0} G$, $\theta=\ga_m$ for some $m\in\aut\d$ such that $m(\al)=\al$. Therefore $b_\al\circ m\circ b_{-\al}\in \aut_0 \d$, and so there exists $\eta\in \t$ such that
\[
m= b_{-\al} \circ \rho_\eta \circ b_{\al}.
\]
Since $m\mapsto \ga_m$ is an isomorphism,
\[
\theta = \ga_{b_{-\al}} \circ \ga_{\rho_\eta} \circ \ga_{b_{\al}}.
\]
It follows by the chain rule that
\begin{align*}
\theta'(s_0) &= X \ga_{\rho_\eta}'{(0,0)}X^{-1}\\
	&\sim X\diag\{\eta,\eta^2\} X\inv
\end{align*}
where $X=\ga_{b_{-\al}}'{(0,0)}$. But $\diag\{\eta,\eta^2\}$
has eigenspaces $\c \oplus 0$ and $0 \oplus \c$ with corresponding eigenvalues $\eta$ and $\eta^2$. Therefore, $\theta '(s_0)$ has eigenspaces $X (\c \oplus 0)$ and $X(0 \oplus \c)$ with corresponding eigenvalues $\eta$ and $\eta^2$.   We have
\[
\ga_{b_{-\al}}(s)=\frac{\left(2\al+(1+|\al|^2)s^1+2\bar\al s^2,s^2+\al s^1+\al^2\right)}{1+\bar\al s^1+\bar\al^2 s^2}.
\]
Hence
\[
X=\ga_{b_{-\al}}'(0,0) \sim (1-|\al|^2) \bbm 1 & 2\bar\al \\ \al &1+|\al|^2 \ebm,
\]
and therefore
\[
X\bpm \c \\ 0\epm= \c\bpm 1 \\ \al \epm =T_{s_0}\calr,\quad X\bpm 0\\ \c \epm = \c\bpm 2\bar\al \\ 1+|\al|^2 \epm = T_{s_0} F_{s_0}.
\]
Thus  $T_{s_0}\calr$ and $T_{s_0}F_{s_0}$ are eigenspaces of $\theta'(s_0)$ with corresponding eigenvalues $\eta,\eta^
2$ respectively.
\end{proof}
\begin{proposition}\label{autlem40}
Let $s_0=(2\al,\al^2)$ for some $\al\in\D$ and let $m\in \aut_{\al} \d$.  If $g$ is any proper embedding of $\d$ into $G$ such that $g(\d) =F_{s_0}$ and $g(\al)=s_0$, then
\[
\ga_m\circ g = g\circ m\circ m.
\]
\end{proposition}
\begin{proof}
Note that 
\[
\ga_m\circ R = R\circ m.
\]
This equation implies that
\[
\ga_m'(s_0) R'(\al) =m'(\al) R'(\al),
\]
which is to say that $m'(\al)$ is the eigenvalue of $\ga_m'(s_0)$ corresponding to the eigenspace $T_{s_0} \calr$. Consequently, by Lemma \ref{autlem30},
\begin{align}\label{aut116}
m'(\al)^2 &\text{ is the eigenvalue of  } \ga_m'(s_0) \notag \\
	&\text{ corresponding to the eigenspace } T_{s_0} F_{s_0}.
\end{align}

Since $g$ is a proper embedding and $\ga_m (F_{s_0}) = F_{s_0}$, there exists $b \in \aut{\d}$ such that
\[
\ga_m\circ g = g\circ b.
\]
 As $b(\al)=\al$, this equation implies that
\[
\ga_m'(s_0) g'(\al) =b'(\al) g'(\al),
\]
that is, $b'(\al)$ is the eigenvalue of $\ga_m'(s_0)$ corresponding to the eigenspace $T_{s_0} F_{s_0}$. Therefore, statement \eqref{aut116} implies that
\[
b'(\al) =m'(\al)^2 = (m\circ m)'(\al),
\]
Since $b, \, m\circ m \in \aut \d$, $b(\al) = \al = (m \circ m)(\al)$, and $b'(\al) = (m\circ m)'(z_0)$, it follows that $b=m\circ m$.
\end{proof}
We describe the phenomena described in Propositions \ref{autlem30} and \ref{autlem40} as the {\em synchrony property} of $G$.
\index{synchrony}

\section{Royal manifolds}\label{royal}

Perhaps the most far-reaching feature of the complex geometry of $G$ is the existence of the special variety $\calr$ with the properties described in Proposition \ref{autprop10}.  We formalize these properties in order to characterize $G$ up to isomorphism.

\begin{defin}\label{autdef10}
Let $\Omega$ be a complex manifold.   We say that $D$ is a \emph{royal disc in $\Omega$}
\index{royal!disc}
 if $D$ is a properly embedded analytic disc in $\Omega$ 
and $D$ satisfies the three conditions of Proposition  \rm{\ref{autprop10}}, that is,
\begin{enumerate}
\item[\rm{(1)}] every automorphism of $\Omega$ leaves $D$ invariant,
\item[\rm{(2)}] every automorphism of $\Omega$ is uniquely determined by its values on $D$,
\item[\rm{(3)}] every automorphism of $D$ has an extension to an automorphism of $\Omega$.
\end{enumerate}
A \emph{royal manifold} 
\index{royal!manifold}
is an ordered pair $(\Omega,D)$
\index{$(\Omega,D)$}
 where $\Omega$ is a complex 
manifold and $D$ is a royal disc in $\Omega$.
\end{defin}

The following lemma is straightforward. 

\begin{lem}\label{autlem4}
If $\Omega$ is a complex manifold and $\Lambda:G\to \Omega$ is a biholomorphic map, then $\Lambda(\calr)$ is a royal disc in $\Omega$ and $(\Omega, \Lambda(\calr))$ is a royal manifold.
\end{lem}
The next proposition spells out the analog of formula \eqref{aut44} on a general royal manifold.
\begin{prop}\label{autprop20}
Let $(\Omega,D)$ be a royal manifold. Then $\aut\Omega$ is isomorphic to $\aut \d$. Furthermore, if $d:\d \to \Omega$ is a properly embedded analytic disc such that $d(\d)=D$, then there exists a unique isomorphism $\Theta:\aut \d\to \aut\Omega$ such that
\be\label{aut50}
\Theta(m) \circ d = d\circ m \qquad \mbox{ for all }  m\in \aut \d.
\ee
\end{prop}

The counterpart of the commutative diagram \eqref{commDiagG} is
\be\label{commDiagO}
\begin{array}{rrrcl}
\d & \stackrel{d}{\longrightarrow} &D& \stackrel{\iota_D}{\longrightarrow} & \Omega \\
\vcenter{\llap{$\scriptstyle{m}$}}\Big\downarrow&  &\vcenter{\llap{$\scriptstyle{\Theta(m)|D}$}}\Big\downarrow & & \vcenter{\llap{$\scriptstyle{\Theta(m)}$}}\Big\downarrow\\
\d & \stackrel{d}{\longrightarrow} & D & \stackrel{\iota_D}{\longrightarrow} & \Omega
\end{array}
\ee
where $\iota_D$ is the injection of $D$ into $\Omega$.

\begin{proof}  
Fix a properly embedded analytic disc $d$ such that $d(\d)=D$. For each $\tau \in \aut\Omega$, Condition (i) in Definition \ref{autdef10} implies that there exists a function $\phi_\tau:\d \to \d$ such that
\be\label{aut60}
\tau\circ d(z) = d\circ\phi_\tau (z) \qquad \mbox{ for } z \in \d.
\ee
Clearly, since $\tau$ is an automorphism of $\Omega$ and $d$ is a properly embedded analytic disc, $\phi_\tau \in \aut \d$.

If $\tau_1,\tau_2 \in \aut\Omega$, then for each $z\in \d$ we see using equation \eqref{aut60} that
\begin{align*}
d(\phi_{\tau_2 \circ \tau_1 }(z))&=\tau_2 \circ \tau_1 \circ d(z)\\
&=\tau_2 ( \tau_1 \circ d(z))\\
&=\tau_2(d\circ \phi_{\tau_1}(z))\\
&=d(\phi_{\tau_2}(\phi_{\tau_1}(z)))\\
&=d(\phi_{\tau_2}\circ \phi_{\tau_1}(z)).
\end{align*}
This relation proves that the map $\Psi:\aut\Omega \to \aut \d$ given by
\be\label{aut70}
\Psi( \tau)= \phi_\tau
\ee
is a homomorphism of automorphism groups.

If $\tau_1,\tau_2 \in \aut\Omega$ and $\phi_{\tau_1}(z) = \phi_{\tau_2}(z) $ for all $z \in \d$, then equation  \eqref{aut60} implies that $\tau_1(d(z))=\tau_2(d(z))$ for all $z\in \d$, which is to say that $\tau_1$ and $\tau_2$ agree on $D$.  Hence, by Condition (ii) in Definition \ref{autdef10}, $\tau_1=\tau_2$. This proves that $\Psi$ is injective.

Consider any $b\in \aut \d$. The map
\[
 d(z) \mapsto d(b(z)) \in D \qquad \mbox{ for } z \in \d
\]
is an automorphism of the complex manifold $D$. Condition (iii) in Definition \ref{autdef10} implies that there exists $\tau \in \aut\Omega$ such that $\tau(d(z)) = d\circ b(z)$ for all $z\in \d$. But then
\[
d(b(z)) = \tau(d(z)) = d(\phi_\tau(z))
\]
for all $z \in \d$, so that $\phi_\tau = b$. This proves that $\Psi$ is surjective from $\aut\Omega$ onto $\aut \d$.

We have shown that $\Psi$ is an isomorphism of groups. In particular, the first assertion of Proposition \ref{autprop20}  (that $\aut\Omega$ is isomorphic to $\aut \d$) is proven. To define an isomorphism $\Theta$ satisfying the second assertion of the proposition,  let $\Theta  =\Psi\inv$. Then $\Theta$ is an isomorphism from $\aut \d$ onto $\aut\Omega$, and equation \eqref{aut50} follows from the relation \eqref{aut60}.

To see that $\Theta$ is unique consider $m\in \aut \d$ and observe that if $\Theta_1$ and $\Theta_2$ are isomorphisms satisfying equation \eqref{aut60}, then $\Theta_1(m)(d(z))=\Theta_2(m)(d(z))$ for all $z \in \d$. Since $\Theta_1(m) $ and $\Theta_2(m)$ agree on $D$, Condition (ii) in Definition \ref{autdef10} imples that $\Theta_1(m) = \Theta_2(m)$. As $m$ is arbitrary, $\Theta_1=\Theta_2$.
\end{proof}

In the light of Proposition \ref{autprop20} we adopt the following definition.
\begin{defin}\label{autdef15}
Let $(\Omega,D)$ be a royal manifold. We say that $(d,\Theta)$ is a \emph{concomitant pair for $(\Omega,D)$}
\index{concomitant pair}
\index{$(d,\Theta)$}
 if  $d:\d \to \Omega$ is a proper analytic embedding, $d(\d)=D$, and $\Theta:\aut \d\to \aut\Omega$ is an isomorphism of groups that satisfies, for all $m\in\aut\d$, 
\[
\Theta(m)\circ d= d\circ m
\]
as in equation \eqref{aut50}.
\end{defin}
  In other words, $(d,\Theta)$ is a concomitant pair for $(\Omega,D)$ if the diagram \eqref{commDiagO} commutes for every $m\in\aut \d$.
\begin{remark}\label{autrem10} \rm 
Concomitant pairs are essentially unique in the following sense. If $(\Omega,D)$ is a royal manifold and $(d_0,\Theta_0)$ is a concomitant pair for $(\Omega,D)$, then $(d,\Theta)$ is a concomitant pair for $(\Omega,D)$ if and only if there exists $b\in \aut \d$ such that $d=d_0 \circ b$ and $\Theta=\Theta_0 \circ I_b$, where $I_b$ denotes the inner automorphism of $\aut\d$ defined by $I_b (m) = b\circ m \circ b^{-1}$.
\end{remark}

As a companion to Lemma \ref{autlem4} we have the following equally straightforward lemma.
\begin{lem}\label{autlem6}
If $\Omega$ is a complex manifold, $\Lambda:G\to \Omega$ is a biholomorphic map, $d=\Lambda\circ R$ and 
\be\label{aut80}
\Theta(m) = \Lambda\circ\gamma_m\circ \Lambda^{-1},\qquad \mbox{ for } m \in \aut{\d},
\ee
then $(d,\Theta)$ is a concomitant pair for $(\Omega,\La(\calr))$.
\end{lem}
\begin{definition}\label{defconsist}
A concomitant pair $(d,\Theta)$ for a royal manifold $(\Omega,D)$ is {\em consistent with}
\index{consistent} 
a bijective map $\La:G\to\Omega$
if $d=\Lambda\circ R$ and $\Theta(m) = \Lambda\circ\gamma_m\circ \Lambda^{-1}$ for all $m\in\aut\d$.
\end{definition}
\subsection{Regularity properties of royal manifolds}\label{regroyal}

\begin{defin}\label{autdef24}
Let $(\Omega,D)$ be a royal manifold and $(d,\Theta)$ a concomitant pair. We say that $(\Omega,D)$ is a \emph{regular} royal manifold if 
\index{royal!manifold!regular}
\begin{enumerate}[\rm (1)]
\item $\Theta:\aut\d\to\aut\Omega$ is differentiable;
\item for every $\lambda\in \Omega \setminus D$, the stabilizer of $\la$ in $\aut\Omega$ is finite, and
\item for every $\lambda\in \Omega \setminus D$,
$ e_\lambda'(\idd)$ is an invertible real-linear map,
where $e_\la:\aut\d\to \Omega$ is defined by 
\beq\label{defela}
e_\la (m)=\Theta(m)(\la).
\eeq
\end{enumerate}
\end{defin}
\begin{remark} \rm
If the complex manifold $\Omega$ is isomorphic to a bounded taut domain \cite{jp},
then $\Theta$ is automatically differentiable -- indeed, by a theorem of H. Cartan \cite{cartan}, real analytic.
\end{remark}
\index{Cartan!H.}
$e_\la'(\idd)$ is a real-linear map between real tangent spaces,
\[
e_\la'(\idd): \laut \to T_\la\Omega.
\] 
 
Conditions (1) to (3) are certainly necessary for $\Omega$ to be biholomorphic to $G$.  They do not depend on the choice of concomitant pair for $(\Omega,D)$.

The following statement is simple to prove.
\begin{proposition}\label{regnec}
If $\Omega$ and $\La$ are as in Lemma {\rm \ref{autlem6}} then $(\Omega, \La(\calr))$ is a regular royal manifold.
\end{proposition}

There is an analog of Proposition \ref{autrem15} for $G$.
\begin{prop}\label{localhomeo}
If $(\Omega,D)$ is a regular royal manifold then, for any $\la\in\Omega\setminus D$, the map $e_\la:\aut\d\to\orb_\Omega(\la)$ is a local homeomorphism and an $N$-to-one covering map, where $N$ is the order of the stabilizer group of $\la$ in $\aut\Omega$.
\end{prop}
\begin{proof}
Let $(d,\Theta)$ be a concomitant pair for $(\Omega,D)$ and
let $H$ be the stabilizer of $\la$ in $\aut\Omega$.  By condition (2) in Definition \ref{autdef24}, $H$ is a finite subgroup of $\aut\Omega$.
For any $m_1,m_2 \in\aut\d$,
\[
e_\la(m_1)=e_\la(m_2) \Leftrightarrow m_2\inv\circ m_1 \in \Theta\inv(H).
\]
Since $\Theta$ is bijective, $|\Theta\inv(H)|=N$.
It follows that $e_\la$ is an $N$-to-one map.

To prove that $e_\la$ is a local homeomorphism, consider any point $e_\la(\beta)$ of $\orb_\Omega(\la)$, where $\beta\in\aut\d$.
Choose a neighborhood $U$ of $\idd$ such that
\[
\{m_2\inv\circ m_1:m_1,m_2 \in U\} \cap \Theta\inv(H) = \{\idd\}.
\]
Let $V$ be a compact neighborhood of $\idd$ contained in $U$.
Then $\beta\circ V$ is a compact neighborhood of $\beta$ on which $e_\la$ is injective, and so $e_\la|\beta\circ V$ is a homeomorphism
onto its range.  Thus $e_\la$ is a local homeomorphism.
\end{proof}
\begin{remark}\label{autrem40} \rm
The natural analog of equation \eqref{aut100}, to wit, 
\[
\rank_\r e_\lambda'(\idd) =2 \quad \mbox{  if }\lambda\in D,
\]
 is not required in  Definition \ref{autdef24}, since the condition holds automatically, as is clear from Proposition \ref{tangD} below.
\end{remark}

\begin{prop}\label{autprop35}
Let $(\Omega,D)$ be a regular royal manifold with concomitant pair $(d,\Theta)$. 
\begin{enumerate}[\rm (1)]
\item  If $\lambda\in D$ then $\orb_\Omega (\lambda)$ is a one-dimensional complex manifold properly embedded in $\Omega$.  
\item  If $\lambda \in \Omega\setminus  D$, then $\orb_\Omega (\lambda)$ is a three-dimensional real manifold properly embedded in $\Omega$.
\end{enumerate}
In either case,
\beq\label{Tlaranela}
\ran e_\la'(\idd)= T_\la\orb_\Omega(\la).
\eeq
\end{prop}

\begin{proof}
(1)  Let $\la\in D$.  By conditions (1) and (3) in Definition \ref{autdef10}, $\orb_\Omega(\la) = D$, which is by hypothesis a properly embedded analytic disc in $\Omega$ and therefore a one-dimensional complex manifold.

(2)   The proof that $\orb_\Omega(\la)$ is a $3$-dimensional real  manifold for any $\la\in\Omega\setminus D$
is almost identical to the proof of the corresponding statement for $G$ in  Theorem \ref{autprop30a}, and so we omit it. 
\end{proof}

For a domain $U$ in $\c^n$, when necessary we shall write $U_r$ for $U$ considered as a $2n$-dimensional real manifold and $U_c$ for $U$ as a complex manifold.  For $p\in U$ the spaces $T_pU_r, \, T_pU_c$ are respectively the real and complex tangent spaces to $U$ at $p$.  We regard elements of $T_pU_r$ as point derivations at $p$ on the algebra $C^1_p(U)$ of germs at $p$ of real-valued $C^1$ functions on $U_r$.  Elements of $T_pU_c$ are point derivations at $p$ on the algebra $\calo_p(U)$ of germs at $p$ of holomorphic functions on $U_c$. We express the action of a point derivation $\de$ on a germ $g$ of the appropriate type by the notation $\ip{g}{\de}$.

The complexification $(T_pM)_\c$ of the real tangent space at $p$ to a real manifold $M$ is the complex vector space comprising the point derivations at $p$ on the complex algebra  $C^1_p(M,\c)$ of germs at $p$ of {\em complex}-valued $C^1$ functions on $M$.
If $\de\in(T_pM)_\c$ then the functional $\re\de$ on $C^1_p(M)$ defined by
\[
\ip{g}{\re\de}= \re \ip{g+i0}{\de}
\]
is a point derivation, that is, a member of $T_pM$.  We also define $\im\de\in T_pM$ to be $-\re(i\de)$. 
In the reverse direction, for a tangent vector $\de\in T_pM$ we denote by $\de_\c$ the complexification of $\de$, so that, for any complex-valued $C^1$ function $h$ in a neighborhood of $p$,
\beq\label{deC}
\ip{h}{\de_\c}= \ip{\re h}{\de}+i\ip{\im h}{\de}.
\eeq
Then, for $\de\in (T_pM)_\c$, the relation $\de= (\re\de)_\c+i(\im\de)_\c$ holds.  Note that, for $\de\in T_pM$, we have $\re (\de_\c) = \de$.  

Furthermore, since every holomorphic function on $U_c$ is a $\c$-valued $C^1$ function on $U_r$, every tangent vector $\de\in  (T_pU_r)_\c$ determines by restriction an element $\de|\calo$ of $T_pU_c$. 

We can summarize the various tangent spaces and their inclusions in the diagram
\beq\label{tgntspcs}
\begin{array}{ccccc}
C^1_p(U_r) & \hookrightarrow & C^1_p(U_r,\c) & \hookleftarrow & \calo_p(U_c) \\
    & & & & \\
T_pU_r & \stackrel{ \stackrel{\cdot_\c}{\rightleftharpoons}}{\scriptstyle{\re,\, \im}}  &  (T_pU_r)_\c  & \stackrel{\cdot|\calo}{\rightarrow} & T_pU_c
\end{array}
\eeq
The vector spaces in the bottom row are respectively real of dimension $2n$, complex of dimension $2n$ and complex of dimension $n$.  The composition of $\cdot_\c$ and $\cdot|\calo$ is a natural real-linear map
\[
\kappa: T_pU_r \to T_pU_c, \quad \mbox{ where } \kappa\de = \de_\c|\calo_p(U_c).
\]
For $\de\in T_pU_r$, the complex tangent vector $\kappa\de$ satisfies, for $g\in\calo_p(U_c)$,
\begin{align}
\ip{g}{\kappa\de} &= \ip{g}{\de_\c} \notag \\
	&= \ip{\re g}{\de}+i\ip{\im g}{\de}, \label{kapaction}
\end{align}
the last line by equation \eqref{deC}.
In terms of the traditional co-ordinates $z^j=x^j+iy^j$ in a neighborhood of $p$, 
\[
\kappa\left(\left(\frac{\partial}{\partial x^j}\right)_p\right)=\left(\frac{\partial}{\partial z^j}\right)_p.
\]
Therefore $\kappa$ is surjective, and since both domain and codomain have real dimension $2n$, it follows that
$\kappa$ is a real linear isomorphism.

  For $\la\in D$ the orbit $\orb_\Omega(\la)$ is the royal disc $D$, which is a properly embedded analytic disc under the complex structure induced by $\Omega$.  Let the evaluation map $e_\la$ be as in Definition \ref{autdef24}, so that $e_\la(m)=\Theta(m)(\la)$. The derivative $e'_\la(\idd)$ is then a real-linear map from $\laut$ to the real tangent space $T_\la D_r$, and so if $\kappa:T_\la D_r\to T_\la D_c$ is the natural embedding of real and complex tangent spaces, then 
\beq\label{kapeprime}
\kappa  e_\la'(\idd):\laut \to T_\la D_c
\eeq
 is a real-linear map from a  $3$-dimensional real space to a $1$-dimensional complex space. 
In fact this map is surjective.
\begin{prop}\label{tangD}
Let $(\Omega, D)$ be a regular royal manifold and let $\la\in D$.  Then
\be\label{aut112}
\ran  \kappa e_\lambda'(\idd)=T_\lambda D_c.
\ee
\end{prop}
\begin{proof}
Let $(d,\Theta)$ be a concomitant pair for $(\Omega, D)$.

  Since $\la\in D$, there exists $z_0\in\d$ such that  $\la=d(z_0)$. Consider a tangent vector $\de$ to $\aut\d$ at $\idd$.  We shall calculate $e_\la'(\idd)\de$.  For any germ $g$ of  real-valued $C^1$ functions on $\Omega$ at $\la$,
\[
\ip{g}{e_\la'(\idd)\de}= \ip{g\circ e_\la}{\de}.
\]
For $m\in\aut\d$,
\begin{align*}
g\circ e_\la(m) &= g\circ \Theta(m)\circ d(z_0) \\
	&=g\circ d\circ m(z_0).
\end{align*}
Hence
\[
\ip{g}{e_\la'(\idd)\de}= \ip{g\circ d\circ m(z_0)}{\de},
\]
where $g\circ d\circ m(z_0)$ is understood as a real-valued function of $m$, with $z_0$ fixed.  Recall the local co-ordinates $r,\al$ for $\aut\d$ introduced in equation \eqref{aut90}.  Here we shall write $\al=\xi+i\eta$, with $\xi,\eta\in\r$ and shall use the local co-ordinates $r,\xi,\eta$ for $\aut\d$.    Note that $\idd\in\aut\d$ corresponds to the local co-ordinates $r=\xi=\eta=0$.   By an elementary calculation,
\begin{align*}
\ip{g\circ d(m_{r,\al}(z_0))}{\left(\frac{\partial}{\partial r}\right)_{\idd}}&=(g\circ d)'(z_0) iz_0, \\
\ip{g\circ d(m_{r,\al}(z_0))}{\left(\frac{\partial}{\partial \xi}\right)_{\idd}}&=(g\circ d)'(z_0)(z_0^2-1),\\
\ip{g\circ d(m_{r,\al}(z_0))}{\left(\frac{\partial}{\partial \eta}\right)_{\idd}}&=(g\circ d)'(z_0)(-i)(z_0^2+1).
\end{align*}
Here $(g\circ d)'(z_0)$ is a real linear functional on $T_{z_0}\d_r$.
   If 
\[
\de= \left(\de_1\frac{\partial}{\partial r} + \de_2 \frac{\partial}{\partial \xi} + \de_3 \frac{\partial}{\partial \eta}\right)_{\idd}
\]
for some real $\de_1,\de_2,\de_3$, then
\be\label{3factors}
\ip{g}{e_\la'(\idd)\de}=(g\circ d)'(z_0) \left(\de_1 iz_0+\de_2(z_0^2-1)-i\de_3(z_0^2+1) \right).
\ee
Now we calculate $\kappa e_\la'(\idd) \in T_\la D_c$.  To this end consider any $h\in\calo_\la(D)$.   By equation \eqref{kapaction},
\begin{align*}
\ip{h}{\kappa e_\la'(\idd)\de} &= \ip{h}{(e_\la'(\idd)\de)_\c} \\
	&= \ip{\re h}{e_\la'(\idd)\de}+i \ip{\im h}{e_\la'(\idd)\de}.
\end{align*}
Thus, by equation \eqref{3factors},
\begin{align*}
\ip{h}{\kappa e_\la'(\idd)\de} &= \ip{\re h}{e_\la'(\idd)\de}+i \ip{\im h}{e_\la'(\idd)\de} \\
	&=\left( (\re h\circ d)'(z_0) +i (\im h\circ d)'(z_0)\right)\left(\de_1 iz_0+\de_2(z_0^2-1)-i\de_3(z_0^2+1) \right) \\
	&= ( h\circ d)'(z_0) \left(\de_1 iz_0+\de_2(z_0^2-1)-i\de_3(z_0^2+1) \right).
\end{align*}
Since $d$ is only determined up to composition with an automorphism of $\d$, no generality is lost by the assumption that $z_0=0$.
Hence
\[
\ip{h}{\kappa e_\la'(\idd)\de} = -(\de_2+i\de_3) (h\circ k)'(0).
\]
On the other hand,
\[
\ip{h}{d'(0)\left(\frac{d}{dz}\right)_0}=\ip{h\circ d}{\left(\frac{d}{dz}\right)_0}=(h\circ d)'(0),
\]
and therefore
\[
\kappa  e_\la'(\idd)\de=-(\de_2+i\de_3) d'(0)\left(\frac{d}{dz}\right)_0.
\]
Thus
\[
\ran \kappa  e_\la'(\idd) =\c d'(0) \left(\frac{d}{dz}\right)_0 = T_\la D_c.
\]
\end{proof}

\subsection{Flat fibrations over royal discs}\label{flatfibO}
In this subsection we shall formalize the consequences for isomorphs of $G$ of the flat fibration of $G$ described in Subsection \ref{flatfibG}.
\begin{defin}\label{autdef30}
Let $(\Omega,D)$ be a royal manifold. If $\cale =\set{E_\lambda}_{ \lambda\in D}$ is a family of subsets of $\Omega$ indexed by $D$, then we say that $\cale$ is a \emph{flat fibration of $\Omega$ over $D$}
\index{flat!fibration}
 if
\begin{enumerate}[\rm (1)]
\item for each $\lambda \in D$, $E_\lambda$ is a properly embedded analytic disc in $\Omega$ such that $E_\la\cap D=\{\la\}$;
\item $\cale$ is a partition of $\Omega$, and  
\item if $\theta \in \aut\Omega$ and $\lambda\in D$, then $\theta (E_\lambda) = E_{\theta(\lambda)}$.
\end{enumerate}
We say that $(\Omega,D,\cale)$ is a \emph{flatly fibered royal manifold}
\index{royal!manifold!flatly fibered}
\index{$(\Omega,D,\cale)$}
 if $(\Omega,D)$ is a royal manifold and $\cale$ is a flat fibration of $\Omega$ over $D$.
We define the {\em flat direction}
\index{flat!direction}
\index{$\lambda^\flat$}
 $\la^\flat$ at a point $\la$ in $\Omega$ to be the tangent space at $\la$ to $E_\mu$ where $\mu\in D$ and $\la\in E_\mu$.
\end{defin}

Clearly, if $(\Omega,D,\cale)$  is a flatly fibered royal manifold then $\Omega$ has complex dimension $2$.

Note that if $(\Omega,D)$ happens to be $(G,\calr)$ then the definition of the flat direction is consistent with that given earlier in Definition \ref{flatdir}.

\begin{lem}\label{autlem24}
Let $\Omega$ be a complex manifold, let  $\Lambda:G \to \Omega$ be a biholomorphic map, let $D=\Lambda(\calr)$ and let
\[
E_{\Lambda(s)} = \Lambda(F_s) \quad \mbox{ for every }s\in\calr,
\]
where $\{F_s:s\in\calr\}$ is the flat fibration of $G$.
Then
\[
\cale=\{E_{\Lambda(s)}: s\in \calr\}
\]
is a flat fibration of the royal manifold $(\Omega,D)$ over $D$, and $(\Omega,D,\cale)$ is a flatly fibered royal manifold.
\end{lem}
\begin{proof}
By Lemma \ref{autlem4}, $(\Omega,D)$ is a royal manifold.  Since $\Lambda$ is a bijection from $\calr$ to $D$, we may write
$\cale=\{E_\la:\la\in D\}$.  Properties (1) and (2) of Definition \ref{autdef30} for the sets $E_\la$ follow from the corresponding properties of the sets $F_s$ for $G$.  If $\theta\in\aut\Omega$ then $\Lambda\inv\circ\theta\circ\Lambda\in\aut G$.  Consider any $s\in\calr$ and $\la=\Lambda(s)\in D$.  We have 
\[
\theta(E_\la)=\theta\circ\Lambda(F_s)= \Lambda\circ(\Lambda\inv\circ\theta\circ\Lambda)(F_s)=\Lambda(F_{\Lambda\inv\circ\theta\circ\Lambda(s)}),
\]
the last step by virtue of property (3) for the flat geodesics as a flat fibration of $(G,\calr)$.  Write $\tilde s=\Lambda\inv\circ \theta(\la)$.
Now
\[
\Lambda(F_{\Lambda\inv\circ\theta\circ\Lambda(s)})=\Lambda(F_{\Lambda\inv\circ\theta(\la)})=\Lambda(F_{\tilde s})=E_{\Lambda(\tilde s)}=E_{\theta(\la)}.
\]
Hence
\[
\theta(E_\la)=E_{\theta(\la)}\quad \mbox{ for all } \la\in D.
\]
Thus the partition $\cale$ has the property (3) of Definition \ref{autdef30}, and so $(\Omega,D,\cale$) is a flatly fibered royal manifold.
\end{proof}

 \subsection{Synchrony in $\Omega$}\label{synchr}

Lemma \ref{autlem40} suggests the following definition concerning the action of $\aut\Omega$ on a flat fibration.
\begin{defin}\label{autdef31}
Let $(\Omega,D,\cale)$ be a flatly fibered royal manifold with concomitant pair $(d,\Theta)$, let $\la_0 \in D$ and let $\la_0=d(z_0)$ for some $z_0\in\D$.


We say that $(\Omega,D,\cale)$ is {\em synchronous at $\la_0$}
\index{synchronous}
 if, for some properly embedded analytic disc
$f:\d\to \Omega$ such that $f(z_0)=\la_0$ and $f(\d)=E_{\la_0}$,
\beq\label{synch}
\Theta(m)\circ f =f\circ m \circ m 
\eeq
for all $m\in \aut_{z_0} \d$.
\end{defin}
\begin{remark}\label{autrem50} \rm
If $(\Omega,D,\cale)$ is as in the definition, then the synchrony of $(\Omega,D,\cale)$  at $\la_0$ depends neither on the choice of $(d,\Theta)$ nor the choice of $f$.
\end{remark}
For  suppose $(\Omega,D,\cale)$ is synchronous  at $\la_0$ with respect to the concomitant pair $(d,\Theta)$ and let $(f_1,\Theta_1)$ be a second concomitant pair.  By Remark \ref{autrem10} there exists $b\in\aut\D$ such that $f_1=f\circ b$ and $\Theta_1=\Theta \circ I_b$, where $I_b(m)= b\circ m\circ b\inv$ for $m\in\aut\D$.  Let $f_1=f\circ b, \, z_1=b\inv(z_0)$.  Then $d_1(z_1)=\la_0=f_1(z_1)$.
Consider $m\in \aut_{z_1}\D$ and $\zeta\in\D$.  Note that $I_b(m)\in\aut_{z_0}\D$, and therefore, from equation \eqref{synch} with $z=b(\zeta)$,
\begin{align*}
\Theta\circ I_b(m)(f\circ b(\zeta))&=f\circ I_b(m) \circ I_b(m)\circ b(\zeta)\\
	&=f\circ b\circ m\circ m(\zeta).
\end{align*}
Hence
\[
\Theta_1(m)(f_1(\zeta))=f_1\circ m\circ m(\zeta).
\]
This shows that synchrony  at $\la_0$ does not depend on the
choice of concomitant pair.

Nor does it depend on the choice of the map $f$.  For suppose that $f_1$ is a second properly embedded analytic disc of $\D$ in $\Omega$ such that $f_1(z_0)=\la_0$ and $f_1(\D)= E_{\la_0}$.  Then there exists $b\in\aut\D$ such that $f_1=f\circ b$ and $b(z_0)=z_0$.  Consider any $m\in\aut_{z_0}\D$ and $\zeta\in\D$.  By equation \eqref{synch},
\[
\Theta(m)\circ f_1(\zeta)=\Theta(m)\left(f\circ b(\zeta)\right)=f\circ m\circ m(b(\zeta)).
\]
Since $\aut_{z_0}\D$ is conjugate in $\aut\D$ to $\aut_0\D$, it is an abelian group.   Hence
\[
\Theta(m)\circ f_1(\zeta)=f\circ b\circ m\circ m(\zeta) = f_1\circ m\circ m(\zeta),
\]
which is the desired relation for $f_1$.

\begin{remark}\label{autrem60} \rm
If $(\Omega,D,\cale)$ is as in the definition, then $(\Omega,D,\cale)$ is synchronous  at a particular $\la_0\in D$ if and only if $(\Omega,D,\cale)$ is synchronous  at $\lambda$ for every $\lambda\in D$. Consequently, it makes sense to say simply that \emph{$(\Omega,D,\cale)$ is synchronous}.
\end{remark}
For  suppose $(\Omega,D,\cale)$ is synchronous  at $\la_0$ with respect to the concomitant pair $(d,\Theta)$ where $\la_0 = d(z_0)$, and let $\la_1 \in D$. 

Suppose $\la_1 = d(z_1)$, for $z_1 \in \D$, and $b(z_0) =z_1$, for $ b\in\aut\D$.
For every $m\in\aut\D$,  we have $\Theta (m) \circ d = d \circ m$. Hence
\[
\Theta(b)(\la_0) = \Theta(b)\circ d(z_0) = d(b (z_0))= d(z_1)=\la_1.
\]
Let 
\[
f_1= \Theta(b) \circ f \circ b^{-1}:  \D \to E_{\la_1}.
\]
Then
$ f_1(z_1) = \la_1 =d(z_1)$.
Consider $m \in \aut_{z_1}\D$ and $\zeta\in\D$.  Then
\begin{align}\label{7.7.1}
\Theta(m)(f_1(\zeta))&=\Theta(m)(\Theta(b)(f\circ b^{-1}(\zeta))) \nn\\
	&=\Theta(m \circ b)(f\circ b^{-1}(\zeta)) \nn\\
 &=\Theta(b)\Theta(b^{-1}\circ m \circ b)(f( b^{-1}(\zeta))).
\end{align}
Since $m$ fixes $z_1= b(z_0)$, $b^{-1}\circ m \circ b$ fixes $z_0$.
By assumption,  $(\Omega,D,\cale)$ is synchronous  at $\la_0$ with respect to the concomitant pair $(d,\Theta)$. Hence
\begin{align}\label{7.7.2}
\Theta(b^{-1}\circ m \circ b)(f( b^{-1}(\zeta)))&=f \circ (b^{-1}\circ m \circ b) \circ (b^{-1}\circ m \circ b) (b^{-1}(\zeta)) \nn\\
	&=f \circ b^{-1}\circ m \circ m (\zeta).
\end{align}
Therefore, by equations \eqref{7.7.1} and \eqref{7.7.2},
\begin{align*}
\Theta(m)(f_1(\zeta))&=
 \Theta(b)\circ f \circ b^{-1}\circ m \circ m (\zeta)\\
&=f_1\circ m\circ m(\zeta).
\end{align*}
Thus $(\Omega,D,\cale)$ is synchronous  at $\la_1$ with respect to the concomitant pair $(d,\Theta)$.

In view of Remarks \ref{autrem50} and \ref{autrem60}, the following statement follows easily from Lemma \ref{autlem40}.
\begin{lem}\label{autlem44}
If $(\Omega,D,\cale)$ is as in Lemma {\rm \ref{autlem24}}, then $(\Omega,D,\cale$) is synchronous.
\end{lem}
\subsection{The sharp direction in $\Omega$}
For a regular royal manifold $(\Omega,D)$ we may define the sharp direction just as we did for $G$ in Definition \ref{defsharp}.
By Proposition \ref{autprop35}, for $\la\in\Omega$ the space $T_\la\orb_\Omega(\la)$ is either a one-dimensional complex subspace (if $\la\in D$) or  a $3$-dimensional real subspace (if $\la\in\Omega\setminus D$) of $T_\la\Omega$.  We may therefore define the space $\la^\sharp$
to be the unique nonzero complex subspace of $T_\la \orb_\Omega(\la)$.  
\index{$\la^\sharp$}
\index{sharp!direction}
In either case
\[
\lambda^\sharp = T_\lambda \orb_\Omega (\lambda) \cap iT_\lambda \orb_\Omega (\lambda).
\]

Covariance of the sharp direction under automorphisms is proved in the same way as Proposition \ref{autprop40}.
\begin{prop}\label{autprop50}
If $\theta\in \aut\Omega$ and $\lambda \in \Omega$ then
\[
  \theta ( \lambda ) ^\sharp=  \theta'(\lambda)\lambda^\sharp.
\]
\end{prop}

\begin{proposition}\label{autlem50}
Let $\Lambda:G \to \Omega$ be a biholomorphic map and let $(d,\Theta)$ be the concomitant pair for $(\Omega,\La(\calr))$ consistent with $\La$.  If $s\in G$ and $\Lambda(s) =\lambda$, then 
\begin{enumerate}
\item[\rm (1)]
$\Lambda'(s)T_s\orb_G(s)=T_\la\orb_\Omega(\la);  $
\item[\rm (2)]
$  \Lambda'(s)s^\sharp= \lambda^\sharp$.
\end{enumerate}
Moreover, if $s\notin\calr$, then $e_s'(\idd)$ is invertible and
\begin{enumerate}
\item[\rm(3)]
$e_\la'(\idd)e_s'(\idd)\inv= \La'(s)\big| T_s\orb_G(s)$;
\item[\rm(4)]
 $e_\la'(\idd)e_s'(\idd)\inv: T_s\orb_G(s) \to T_\la\orb_\Omega(\la)$ is a real linear map whose restriction to $s^\sharp$ is complex linear and maps $s^\sharp$ to $\la^\sharp$.
\end{enumerate}
\end{proposition}

For $s\in\calr$, the real linear map $e_s'(\idd)$ maps the $3$-dimensional space $\laut $ to $T_s\calr$, which is $2$-dimensional, so we cannot form $e_s'(\idd)\inv$.

\begin{proof}
\nin (1) By assumption, $d=\Lambda\circ R$ and $\Theta(m)=\Lambda\circ \ga_m\circ\Lambda\inv \in\aut\Omega$ for every $m\in\aut\D$.  Hence
\begin{align*}
e_\la(m) &= \Theta(m)(\la)\\
	&= \Lambda\circ \ga_m\circ\Lambda\inv(\la) \\
	&=\Lambda\circ \ga_m(s)\\
	&=\Lambda\circ e_s(m).
\end{align*}
That is, $e_\la=\Lambda\circ e_s$.  Hence
\begin{align}
e_\la'(\idd) &= \Lambda'(e_s(\idd))e_s'(\idd) \notag \\
	&= \Lambda'(s)e_s'(\idd). \label{elaLaes}
\end{align}
Therefore
\[
\ran e_\la'(\idd) = \ran \Lambda'(s)e_s'(\idd) = \Lambda'(s) \ran e_s'(\idd),
\]
which is to say (by virtue of equations \eqref{Tsranes} and \eqref {Tlaranela}) that  
\[
T_\la\orb_\Omega(\la)=\Lambda'(s)  T_s\orb_G(s).
\]

\nin (2)
$s^\sharp$ is a nonzero complex subspace of $T_s\orb_G(s)$.  Since $\Lambda'(s)$ is a nonsingular complex linear map,  $\Lambda'(s)s^\sharp$
is a nonzero complex linear subspace of $T_\la\orb_\Omega(\la)$, hence is $\la^\sharp$.

\nin (3)  Consider $s\in G\setminus \calr$.  By Theorem \ref{autprop30a}, the real linear map $e_s'(\idd)$ has full rank between the $3$-dimensional spaces $\laut $ and $T_s\orb_G(s)$, and so is nonsingular.  Hence, $e_\la'(\idd)e_s'(\idd)\inv$ exists and is a real linear map from $T_s\orb_G(s)$ to $T_\la\orb_\Omega(\la)$.  By equation \eqref{elaLaes},
\beq\label{eeLa}
e_\la'(\idd)e_s'(\idd)\inv = \La'(s) \quad \mbox{ on }\quad T_s\orb_G(s).
\eeq

\nin (4)
Since $\La'(s)$ is a complex linear map on $\c^2$, it follows that $e_\la'(\idd)e_s'(\idd)\inv$ is a complex linear map on the complex linear subspace $s^\sharp$ of $\c^2$. By (2) and equation \eqref{eeLa}, $e_\la'(\idd)e_s'(\idd)\inv s^\sharp=\la^\sharp$.

\end{proof}

\subsection{Sharpness of the action of $\aut\Omega$}\label{sharpness}

In this section, for a flatly fibered royal manifold $(\Omega,D,\cale)$,
we shall show that sharp action of $\aut\Omega$, as described in the introduction, is necessary for $\Omega$ to be isomorphic to $G$. In the next subsection we shall show that the condition is also sufficient.
We first define sharpness more formally than in the introduction.
Recall that, for a flatly fibered royal manifold $(\Omega,D,\cale)$, we defined the Poincar\'e parameter $P(\mu)$ for $\mu\in\Omega$  to be the Poincar\'e distance of $\mu$ from $\la$, where $\la\in D$ and $\mu\in E_\la$ and the distance is taken in the disc $E_\la$.  That is, if $f:\d\to \Omega$ is a proper analytic embedding with range $E_\la$ and $f(z_0)=\la, \, f(z)=\mu$, then
\beq\label{defP}
P(\mu) \; \df \;\arctanh\left| \frac{z-z_0}{1-\bar z_0 z}\right|.
\eeq
\index{$P(\mu)$}
\index{Poincar\'e parameter} 
It will be convenient to use also the pseudohyperbolic variant $C(\mu)$ of $P(\mu)$, defined for $\mu\in E_\la$ to be the pseudohyperbolic distance in $E_\la$ from $\mu$ to $\la$.  In other words, if $f\in\Omega(\d)$ has range $E_\la$ and $f(z_0)=\la, \, f(z)=\mu$, then
\beq\label{defC}
C(\mu)  \df \left|\frac{z-z_0}{1-\bar z_0 z}\right|.
\eeq
Thus $P$ and $C$ are related by the equations
\begin{align}
P(\mu)&= \arctanh C(\mu)\notag \\
	&= \half \log \frac{1+C(\mu)}{1-C(\mu)}. \label{PandC}
\end{align}

\begin{remark}\label{invarC}\rm
Observe that $P(\cdot)$ and $C(\cdot)$ are invariant under isomorphisms which preserve foliations.  If $(\Omega_j,D_j,\cale_j)$ is a flatly fibered royal manifold for $j=1,2$, if $\La:\Omega_1\to\Omega_2$ is an isomorphism which maps the leaves of $\cale_1$ to those of $\cale_2$ and if $\mu\in\Omega_1$ then $C(\mu)=C(\La(\mu))$.
\end{remark} 
\begin{definition}\label{newsharp}
Let $(\Omega,D,\cale)$ be a regular flatly fibered royal manifold having a concomitant pair $(d,\Theta)$.  Let $\mu\in\Omega\setminus D$ and let $(U,\psi)$ be a chart in $\Omega$ such that $\mu\in U$.  We say that $\aut\Omega$ acts sharply at $\mu$ with respect to $(d,\Theta)$ if
\beq\label{shar}
\e^{2P(\mu)}\left(\psi( \Theta(B_{it})(\mu)) - \psi(\mu)\right)= i \left(\psi(\Theta(B_t)(\mu)) - \psi(\mu)\right) + o(t)
\eeq
as $t\to 0$ in $\r$.
\end{definition}
\index{sharp!action}
The condition \eqref{shar} states that the tangents $v_{\bf i}$ and $v_{\bf 1} \in \c^2$ at $t=0$ to the curves $\psi(\Theta(B_{{\bf i}t})(\mu))$ and $ \psi(\Theta(B_{t})(\mu))$ in $\psi(U)$ satisfy 
\[
\e^{2P(\mu)}v_{\bf i}={\bf i}v_{\bf 1},
\]
 where $\bf i$ is (temporarily, for this sentence) the imaginary unit.  This property is clearly independent of the chart $\psi$ since the derivative of any transition function at a point is a complex-linear map.

We need to examine how sharpness depends on the concomitant pair $(d,\Theta)$.
\begin{proposition}\label{relconcom}
With the notation of Definition {\rm \ref{newsharp}}, let $(d_1,\Theta_1)$ (for some $b\in\aut\d$) be the concomitant pair
\[
(d\circ b, \Theta\circ I_b) \quad \mbox{ where } I_b(m)= b\circ m\circ b\inv \; \mbox{ for } m\in\aut\d.
\]
 Let $\mu_1=\Theta(b)(\mu) \in \Omega\setminus D$.
Then $\aut\Omega$ acts sharply at $\mu_1$ with respect to $(d_1,\Theta_1)$ if and only if 
$\aut\Omega$ acts sharply at $\mu$ with respect to $(d,\Theta)$.
\end{proposition}

\begin{proof}
We may choose the chart 
\[
\psi_1= \psi\circ\Theta(b)\inv\quad \mbox{ on } \Theta(b) (U)
\]
at $\mu_1$.  For small real $t$,
\begin{align*}
\psi_1(\mu_1)&=\psi\circ\Theta(b)\inv ((\Theta(b)(\mu))\\
	&= \psi(\mu),\\
\psi_1( \Theta_1(B_{it})(\mu_1)) &=\psi\circ\Theta(b)\inv\circ (\Theta(b)\circ\Theta(B_{it})\circ\Theta(b)\inv)((\Theta(b)(\mu)) \\
	&=\psi(\Theta(B_{it})(\mu)).
\end{align*}
These equations, together with the analogous ones with $it$ replaced by $t$ and the fact that $P(\mu_1)=P(\mu)$, imply the statement in the proposition.
\end{proof}
\begin{definition}\label{newsharp2}
Let $(\Omega,D,\cale)$ be a regular flatly fibered royal manifold having a concomitant pair $(d,\Theta)$.  Let $\mu\in\Omega\setminus D$ and let $(U,\psi)$ be a chart in $\Omega$ such that $\mu\in U$.  We say that $\aut\Omega$ {\em acts sharply on} $\Omega$ if $\aut\Omega$ acts sharply with respect to $(d,\Theta)$ at every point of $\Omega\setminus D$.
\end{definition}
\begin{remark}\label{indep} \rm
(1) Propositions \ref{autrem10} and \ref{relconcom} show that the sharpness of the action of $\aut\Omega$ on $\Omega$ does not depend on the choice of the concomitant pair $(d,\Theta)$.\\
(2) With respect to a fixed concomitant pair, for any $\theta\in\aut\Omega$, automorphisms act sharply at $\mu\in \Omega\setminus D$ if and only if they act sharply at $\theta(\mu)$. Since every orbit in $\Omega$ meets every leaf in $\cale$, to conclude that $\aut\Omega$ acts sharply, it is enough to show that, for some $\la\in D$,  automorphisms act sharply at every point of $E_\la\setminus \{\la\}$.
\end{remark}

\begin{prop}\label{autprop60}
Let $\Omega$ be a complex manifold and $\Lambda:G\to \Omega$ be a biholomorphic map.  There exist a royal disc $D$ in $\Omega$, a flat fibration $\cale$ of $\Omega$ over $D$ and a concomitant pair $(d,\Theta)$ such that $(\Omega,D,\cale)$ is a synchronous regular flatly fibered royal manifold, $(d,\Theta)$ is consistent with $\La$  and $\aut\Omega$ acts sharply on $\Omega\setminus D$.
\end{prop}
\begin{proof}
Let $D=\La(\calr)$.  By Lemmas \ref{autlem4} and  \ref{autlem6}, $(\Omega,D)$ is a royal manifold and there is a concomitant pair $(d,\Theta)$ for $(\Omega,D)$ such that 
\beq\label{consis}
\Theta(m)\circ \La=\La\circ\ga_m
\eeq
for all $m\in\aut\d$.  Thus $(d,\Theta)$ is consistent with $\La$.
By Proposition \ref{regnec}, $(\Omega,D)$ is a regular royal manifold.
Let $ \cale$ correspond under $\Lambda$ to the flat fibration of $G$. By Lemma \ref{autlem44},
$(\Omega,D,\cale)$ is a synchronous regular flatly fibered royal manifold.  It remains to show that $\aut\Omega$ acts sharply on $\Omega\setminus D$.

Consider a point $\mu\in\Omega\setminus D$, say $\mu\in E_\la$, for some $\la\in D$.   Let $s=\La\inv(\mu)$.  We may assume (by modifying $\La$ and $\Theta$ and utilising Remark \ref{indep}) that $s$ has the form $(0,p)$ for some $p\in (0,1)$.  
Then $s$ lies in the flat geodesic $F^0$, and so $\La\inv(E_\la)=F^0$.  It follows that $\La\inv(\la)=(0,0)$ and
since isomorphisms preserve the M\"obius distance, $p=C(\mu)$.

Let $(U,\psi)$ be a chart at $\mu$.
For any $\al\in\c$ and all small enough real $t$, in view of equation \eqref{consis},
\begin{align}\label{tgtmu}
\psi\left(\Theta(B_{t\al})(\mu)\right) &=\psi\circ \La\circ \ga_{B_{t\al}} (s)  \notag\\
&= \psi(\mu) + (\psi\circ\La)'(s) \frac{\dd}{\dd t} \ga_{B_{t\al}}(s) \big|_{t=0}+o(t).
\end{align}
By Lemma \ref{vralpha}, with $r=0$ and $s=(0,p)$,
\[
 \frac{\dd}{\dd t} \ga_{B_{t\al}}(s) \big|_{t=0} = -\al\bpm 2\\0\epm -\bar\al \bpm 2p\\ 0 \epm.
\]
Let $A=-2 (\psi\circ\La)'(s)$, so that $A$ is a complex-linear map.
Taking successively $\al=i$ and $\al=1$ in equation \eqref{tgtmu} we obtain
\begin{align*}
\psi\left(\Theta(B_{ti})(\mu)\right) - \psi(\mu) &= i(1-p)A(1,0)+ o(t), \\
\psi\left(\Theta(B_{t})(\mu)\right) - \psi(\mu) &= (1+p)A(1,0)+o(t).
\end{align*}
Hence
\[
(1+p)\left(\psi\left(\Theta(B_{ti})(\mu)\right) - \psi(\mu)\right)=i(1-p)\left(\psi\left(\Theta(B_{t})(\mu)\right) - \psi(\mu)\right) +o(t).
\]
Since $p=C(\mu)$, this is to say that $\aut\Omega$ acts sharply at $\mu$.
\end{proof}
The next statement justifies the terminology of `sharp action'.
\begin{lemma}\label{oldsharp}
Let $(\Omega,D,\cale)$ be a regular flatly fibered royal manifold 
 and suppose that $\aut\Omega$ acts sharply at a point $\mu\in\Omega\setminus D$.
For any $s\in G$ such that $C(s)=C(\mu)$ the map
\beq\label{sharpclinear}
e_\mu'(\idd)e_s'(\idd)\inv: T_s\orb_G(s) \to T_\mu\orb_\Omega(\mu)
\eeq
maps $s^\sharp$ to $\mu^\sharp$ and is a complex-linear map.
\end{lemma}
\begin{proof}
Let 
\[
X=e_\mu'(\idd) e_s'(\idd)\inv.
\]
$X$ is a real-linear map from $T_sG$ to $T_\mu\Omega$.  We must show that 
$Xs^\sharp\subseteq\mu^\sharp$ and that $X$ is complex-linear on $s^\sharp$.

We can assume that $s=(0,p)$ where $0< p<1$.  Clearly
\[
p=C(s)=C(\mu).
\]
By equation \eqref{PandC},
\[
\e^{2P(\mu)}=\frac{1+p}{1-p}.
\]
The sharpness hypothesis, according to Definition \ref{oldsharp}, is
\beq\label{shnow}
\frac{1+p}{1-p}\left(\psi( \Theta(B_{it})(\mu)) - \psi(\mu)\right)= i \left(\psi(\Theta(B_t)(\mu)) - \psi(\mu)\right) + o(t)
\eeq
as $t\to 0$ in $\r$.
  By Proposition \ref{formsharp},
\[
s^\sharp= \c\bpm 1\\0\epm.
\]

We shall use the local co-ordinates $(r,\al)\in (-\pi,\pi)\times \d$ for a neighborhood of $\idd$ in $\aut\d$, as in Lemma \ref{atlas}.
By Lemma \ref{vralpha},
\[
e_s'(\idd)(\ph_1\inv)'(r,\al)=v_{r,\al}(0,p)=\bpm -2\al-2p\bar\al \\ 2irp \epm,
\]
which is in $s^\sharp$ if and only if $r=0$.  Thus
\[
e_s'(\idd)\inv s^\sharp= (\ph_1\inv)'(0\oplus \c).
\]
Moreover, for all $\al\in\c$,
\beq\label{esprinv}
(\ph_1\inv)'(0,\al)=e_s'(\idd)\inv\bpm -2\al-2p\bar\al \\ 0 \epm.
\eeq

Note that $m_{0,\al} = B_\al$ in the notation of equation \eqref{defBal}.
Let $\psi$ be a chart on $\Omega$ at $\mu$.
For any $\al\in\c$, as $t\to 0$ in $\r$, 
\begin{align*}
\psi\left(\Theta(B_{t\al})(\mu)\right)&=\psi(\mu)+\psi'(\mu)e_\mu'(\idd)(\ph_1\inv)'(0,\al) + o(t) \\
	&=\psi(\mu)+\psi'(\mu)X\bpm -2\al-2p\bar\al \\ 0 \epm    + o(t)
\end{align*}

Take in succession $\al=1$ and $\al=i$ and use the real-linearity of $X$ to obtain the relations
\begin{align}
\psi\left(\Theta(B_{t})(\mu)\right)-\psi(\mu)&=-2(1+p)\psi'(\mu)X (1,0)   + o(t),\label{succ1} \\
\psi\left(\Theta(B_{ti})(\mu)\right)-\psi(\mu)&=-2(1-p)\psi'(\mu)X i(1,0) + o(t).\label{succ2}
\end{align}
We have
\begin{align*}
\psi'(\mu)Xi(1,0)&= -\frac{1}{2(1-p)}\left(\psi(\Theta(B_{ti}(\mu))-\psi(\mu)\right)+o(t)\\
	&\hspace*{4.1cm} \quad\mbox{by equation \eqref{succ2}}\\
	&= -\frac{1}{2(1+p)}\left(\psi(\Theta(B_{t}(\mu))-\psi(\mu)\right)+o(t)\\
	&\hspace*{4.1cm} \quad\mbox{by equation 
\eqref{shnow}}\\
	&=i \psi'(\mu)X(1,0)+o(t)\quad\mbox{by equation \eqref{succ1}}.
\end{align*}
Since $\psi'(\mu)$ is an invertible complex-linear map which identifies $T_\mu\Omega$ with $\c^2$, it follows that
\beq\label{XiiX}
X(i,0)=iX(1,0).
\eeq

The vectors $(1,0)$ and $(i,0)$ span $s^\sharp$ over $\r$, and 
\[
\ran X \subseteq\ran e_\mu'(\idd)\subseteq T_\mu\orb_\Omega(\mu).
\]
Equation \eqref{XiiX} now shows both that $\ran X\subset \mu^\sharp$ and that $X$ is complex-linear on $s^\sharp$.
\end{proof}

\subsection{A characterization of $G$}
We have arrived at the main theorem of the paper.
\begin{theorem}\label{autprop70}
A complex manifold $\Omega$ is isomorphic to $G$ if and only if there exist
a royal disc $D$ in $\Omega$ and a flat fibration $\cale$ of $\Omega$ over $D$ such that 
 $(\Omega,D,\cale)$ is a synchronous regular flatly fibered royal manifold 
 and $\aut{\Omega}$ acts sharply on $\Omega$.
\end{theorem}
\index{main theorem}
\begin{figure}[!ht]
\includegraphics{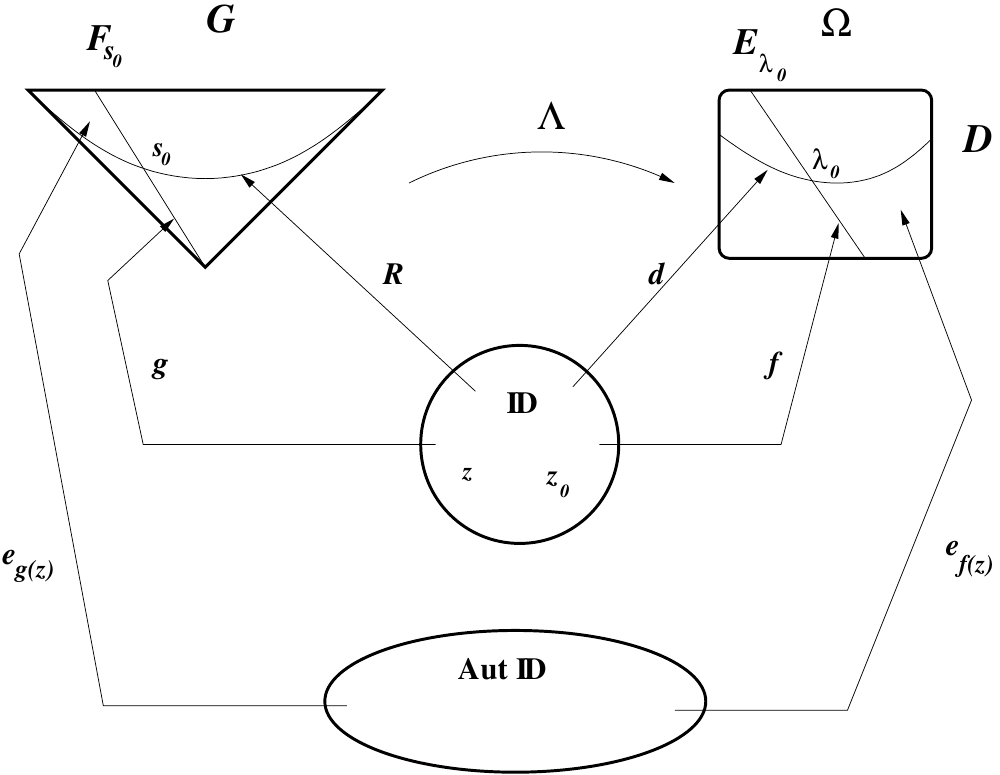}
\caption{The construction of $\La:G\to\Omega$}
\end{figure}
\index{Figure 1}
\begin{proof}
Necessity is Proposition \ref{autprop60}.  We prove sufficiency.
Let  $(d,\Theta)$ be a concomitant pair for $(\Omega,D)$.
Choose $z_0 \in \d$ and let $s_0=R(z_0)$ and $\la_0=d(z_0)$. We shall construct a biholomorphic map $\Lambda:G \to \Omega$ satisfying $\Lambda(s_0) = \la_0$.
Figure 1 represents the construction.

Choose a properly embedded analytic disc $g$ of $\d$ into $G$ satisfying $g(\d) = F_{s_0}$ and $g(z_0) =s_0$.  Choose also a properly embedded analytic disc $f:\d \to \Omega$ such that $f(z_0) = \la_0$ and $f(\d) = E_{\la_0}$. For $s\in G$ we define $\Lambda(s)$ by the following recipe.  

Since each point in $G$ is in a flat geodesic and $\aut  G$ acts transitively on the flat geodesics, we may choose $m\in \aut\d$ such that $\ga_m\inv(s) \in F_{s_0}$ and hence there exists $z\in\d$ such that 

\be\label{getmz}
s=\gamma_m\circ g(z).
\ee
Let
\be\label{aut130}
\Lambda(s)=\Theta(m)\circ f(z)
\ee
Certainly $\La(s)\in\Omega$.
To see that this recipe does define $\Lambda$ as a map from $G$ to $\Omega$, consider $z_1,z_2 \in \d$ and $m_1,m_2 \in \aut \d$ such that
\be\label{aut140}
\gamma_{m_1}\circ g(z_1) = \gamma_{m_2}\circ g(z_2).
\ee
We wish to show that
\be\label{aut150}
\Theta(m_1)\circ f(z_1)=\Theta(m_2)\circ f(z_2).
\ee
Note first that equation \eqref{aut140} implies that if $ m=m_2^{-1}\circ m_1 $, then 
\[
\gamma_m\circ g(z_1)=g(z_2).
\]
 Since $g(z_1), g(z_2) \in F_{s_0}$, Lemma \ref{autlem20} implies that $m\in \aut_{z_0}\d$. Consequently, by Lemma \ref{autlem40}
\[
g\circ m\circ m(z_1) = \gamma_m\circ g(z_1) =g(z_2),
\]
which implies that
\[
z_2 = m \circ m (z_1).
\]
By hypothesis, $(\Omega,D,\cale)$ is synchronous.  According to Definition \ref{autdef31}, it means (since $m\in \aut_{z_0}\d$) that
\[
\Theta(m)\circ f = f\circ m\circ m.
\]
Hence
\begin{align*}
\Theta(m)\circ f(z_1) & = f\circ m\circ m(z_1)\\
	&=f(z_2).
\end{align*}
Therefore equation \eqref{aut150} is true, and so $\La(s)$ is unambiguously defined.

On taking $m=\idd$ in equation \eqref{aut130} we have
\be\label{aut170}
\Lambda \circ g(z) = f(z)
\ee
for all $z\in \d$. In particular, $\Lambda(s_0) =\la_0$.

Consider any $\ups\in\aut\d$.  Since 
\[
\ga_\ups(s)=\ga_\ups\circ \ga_m\circ g(z) = \ga_{\ups\circ m}\circ g(z),
\]
by the definition \eqref{aut130} of $\La$,
\begin{align*}
\La\circ \ga_\ups(s) &= \Theta(\ups\circ m)\circ f(z) \\
	&= \Theta(\ups)\circ\Theta(m)\circ f(z) \\
	&=\Theta(\ups) \circ\La(s).
\end{align*}
Thus
\beq\label{aut160}
\La\circ \ga_\ups= \Theta(\ups)\circ \La \quad \mbox{ for all }\ups\in\aut\d.
\eeq

Now fix a general point $R(z)$ on the royal geodesic $\calr$. If  $m\in\aut\d$ is such that $m(z_0)=z$, then 
\begin{align*}
\Lambda\circ R(z)&= \Lambda\circ R\circ m(z_0) &  \\ 
	&=\Lambda\circ \gamma_m \circ R(z_0) & \mbox{ by equation \eqref{aut44}}.
\end{align*}
 By equation \eqref{aut160}  and the fact that $s_0=R(z_0)$,
\begin{align*}
\Lambda\circ R(z) &=\Theta(m)\circ \Lambda(s_0) &  \\
	&=\Theta(m)( \la_0) &  \\
	&=\Theta(m)\circ d(z_0) &   \\
	&=d\circ m(z_0) & \mbox{ by equation \eqref{aut50}} \\
	&=d(z). & 
\end{align*}
Thus
\be\label{aut180}
\Lambda\circ R = d.
\ee

Now fix $z_1 \in \d$ and choose $m$ such that $m(z_0) = z_1$. Since $\gamma_m\circ R(z_0) = R(z_1)$, Lemma \ref{autlem20} implies that
\[
\gamma_m(F_{R(z_0)})=F_{R(z_1)},
\]
and since $\Theta(m) \circ d(z_0) =d(z_1)$, Condition (3) in Definition \ref{autdef30} implies that
\[
\Theta(m)(E_{d(z_0)})=E_{d(z_1)}.
\]
Therefore
\begin{align*}
\Lambda(F_{R(z_1)}) &= \Lambda\circ \gamma_m(F_{R(z_0)})\\
	&=\Theta (m) \circ \Lambda (F_{R(z_0)})\\
	&=\Theta(m)( E_{d(z_0)})\\
	&=E_{d(z_1)}.
\end{align*}
To summarize, we have shown that if $\calf$ denotes the partition of $G$ in Lemma \ref{autlem10} and $\cale$ denotes the partition of $\Omega$ in Definition \ref{autdef30}, then $\Lambda$ induces a map $\La^\sim: \calf\to\cale$ given by
\[
\Lambda^\sim(F_{R(z)}) = E_{d(z)}.
\]
Furthermore, as the map $R(z) \mapsto d(z)$ from $\calr$ to $D$ 
is a bijection, so also is $\Lambda^\sim$.

Consider any point $\la\in E_{d(z_1)}$.  Then $\Theta(m\inv)(\la)\in E_{\la_0}$, and so $\Theta(m\inv)(\la)=f(z)$ for some $z\in\d$.
Hence $\la=\Theta(m)\circ f(z)$.  By equations \eqref{getmz} and \eqref{aut130}, $\la=\La\circ \la_m\circ g(z)$.  Thus $\la\in\La(G)$, and so $\La$ is surjective.

Suppose $s_1,\, s_2\in G$ satisfy $\La(s_1)=\La(s_2)$.
Since $\Lambda^\sim$ is a bijection, it follows that $s_1,\, s_2$ lie in the same flat geodesic in $G$, say in $F_{R(z_1)}$. Let $m\in\aut\d$  be such that $m(z_0)=z_1$.  We have, for $j=1,2$,
\[
\ga_m\inv(s_j)\in F_{R\circ m\inv(z_1)}=F_{R(z_0)}=F_{s_0}.
\]
Hence  $\ga_m\inv(s_j)=g(\zeta_j)$ for some $\zeta_1,\zeta_2\in\d$.  By equation \eqref{aut130},
\[
\La(s_j)=\Theta(m)\circ f(\zeta_j).
\]
Hence $\Theta(m)\circ f(\zeta_1)=\Theta(m)\circ f(\zeta_2)$, and therefore $\zeta_1=\zeta_2$.  Thus 
\[
s_1=\ga_m\circ g(\zeta_1)=\ga_m\circ g(\zeta_2)=s_2.
\]
We have shown that $\La:G\to\Omega$ is bijective.  Moreover, we can observe that
\be\label{aut190}
\Lambda| F_{R(z_1)} = \Theta (m) \circ g \circ f^{-1} \circ \gamma_m^{-1}| F_{R(z_1)}.
\ee

There remains to prove that $\Lambda$ and $\Lambda^{-1}$ are holomorphic. 

We shall first show that $\La$ is smooth as a mapping between real manifolds by giving a formula for $\La$ which is clearly differentiable.
The assumption that $z_0=0$, $g(z) =(0,z)$ and so $s_0=(0,0)$ loses no generality.  It implies that $F_{s_0}=\{(0,z):z\in\D\}$.

Consider a point
\[
s=(\zeta+\eta,\zeta\eta) \in G
\]
for some $\zeta,\eta\in\D$.  To evaluate $\La(s)$ we shall choose an automorphism $m$ of $\D$ satisfying $m'(0)>0$  such that
$\ga_{m}\inv(s) \in F_{s_0}$.   To see that this is possible take $m=B_\al$ for some $\al\in\D$.  Then $m'(0)>0$.  We require
$\ga_{B_\al}\inv(s) \in F_{s_0}$, 
which is to say that
\[
B_{-\al}(\zeta)+B_{-\al}(\eta)=0.
\]
Expressing this relation in terms of the components $s^1, s^2$ of $s$, we must find $\al= \al(s) \in\D$ such that
\[
s^1=-\frac{2\al}{1+|\al|^2}-\frac{2\bar\al}{1+|\al|^2} s^2.
\]
Compare this expression with that of the flat co-ordinates for $s$ given in equations \eqref{flatbeta} and \eqref{forbeta}:
\[
s^1=\beta + \bar\beta s^2,
\]
where
\[
\beta=\beta(s)=\frac{s^1-\bar{s^1}s^2}{1-|s^2|^2}.
\]
One sees that it suffices to choose $\al(s)$ such that
\[
\beta(s)= -\frac{2\al(s)}{1+|\al(s)|^2}.
\]
A suitable choice of $\al(s)$ is
\[
\al(s)= \frac{-\beta(s)}{1+\sqrt{1-|\beta(s)|^2}}
\]
as may readily be checked.  Clearly  $\beta, \al \in\D$ and both $\beta$ and $\al$ are real-analytic functions of  $s$.
Moreover
\begin{align*}
L(s)\df \ga_{B_{\al}}\inv(s)&=(0, B_{-\al}(\zeta)B_{-\al}(\eta) )\\
	&=\left(0,\frac{s^2+\al(s)s^1+\al(s)^2}{1+\overline{\al(s)}s^1+ \overline{\al(s)}^2s^2}\right),
\end{align*}
which is also real-analytic in $s$.  By the definition of $\La$,
\[
\La(s)= \Theta(B_{\al(s)})\circ f \circ g\inv \circ L(s).
\]
The map $s\mapsto B_{\al(s)}$ is real-analytic from $G$ to $\aut\D$.  
Since the action of $\aut \D$ on $\Omega$ is differentiable, by the regularity assumption on the royal manifold $(\Omega,D)$, we conclude that $\La: G\to \Omega$ is differentiable.

Consider $s\in G$ and suppose that $s\in F_{R(z_1)}$.
Let $X=\Lambda'(s)$ viewed as a real-linear mapping from $T_sG$ to $T_{\Lambda(s)}\Omega$.

  Recall from Definition \ref{flatdir} that $s^\flat$ denotes the flat direction at $s$.   
Equation \eqref{aut190} implies that
\be\label{aut200}
 X(s^\flat) = \Lambda (s)^\flat \text{ and } X|s^\flat \text{ is complex linear}.
 \ee

By equations \eqref{getmz} and \eqref{aut130}, for all $z\in\d$ and $m\in\aut\d$,
\[
\La\circ \ga_m\circ g(z)=\Theta(m)\circ f(z).
\]
In view of the definitions \eqref{defes} and \eqref{defela}, this equation can be written
\[
\La\circ e_{g(z)}=e_{f(z)} :\aut\d\to\Omega.
\]
On differentiating at $\idd$ we obtain
\[
\La'\circ e_{g(z)}(\idd) e_{g(z)}'(\idd)=e_{f(z)}'(\idd):\laut \to T_{f(z)}\orb_\Omega(f(z)).
\]
  For $z\neq z_0$, the point $g(z)\notin \calr$, and therefore, by Proposition \ref{autlem50}, $e'_{g(z)}(\idd)$ is invertible, and so
\[
\La'\circ g(z) =e_{f(z)}'(\idd)e_{g(z)}'(\idd)\inv.
\]
By the hypothesis, $\aut\Omega$ acts sharply on $\Omega$. By Lemma \ref{oldsharp}, it follows that $\La'\circ g(z)$ maps $g(z)^\sharp$ into $f(z)^\sharp$ and is complex-linear on $g(z)^\sharp$  whenever $g(z) \notin \calr$.

Recalling that $X:T_sG\to T_{\Lambda(s)}\Omega$ is real-linear and that (by Proposition \ref{flneqsh}) $s^\flat$ and $s^\sharp$ are linearly independent, we infer from equation \eqref{aut200} that $X=\La'(s)$ is complex linear  for all $s\in G\setminus\calr$.   Therefore  $\Lambda$ is analytic  on $G\setminus\calr$.  
The restriction of $\La$ to any co-ordinate plane $P_\zeta\df\{s\in G: s^1=\zeta\}$, for $|\zeta|<2$, is analytic in $s^2$ except possibly at the sole point $(\zeta, \tfrac 14\zeta^2)$ of $P_\zeta\cap\calr$ and is continuous on $P_\zeta$.  Hence $\La|P_\zeta$ is analytic in $s^2$.  Likewise
the restriction of $\La$ to any of the orthogonal co-ordinate planes is analytic in $s^1$.  Thus $\La$ is analytic on $G$.  Every bijective holomorphic map between domains  has a holomorphic inverse (for example, \cite[Chapter 10, Exercise 37]{kr}).  It follows easily that a bijective holomorphic map between a domain and a complex manifold has a holomorphic inverse.
\end{proof}

\section{A characterization of $G$ via flat co-ordinates} \label{charactflat}

Recall from Subsection \ref{flatfibG} that $G$ is foliated by the sets
\[
F^\beta= \{(\beta+ \bar\beta z,z):z\in\d\}
\]
 for $\beta\in\d$ \cite[Theorem 2.1]{ay2004}. 
Thus the map
$\eta:\d^2 \to G$ defined by the formula
\be\label{30}
\eta(\beta,z) =(\beta +\bar \beta z, z),\qquad \beta,z \in \d,
\ee
is a homeomorphism of $\d^2$ onto $G$.

We will call $\beta, z$ the {\em flat co-ordinates} for points of $G$.
\index{flat co-ordinates}
In this section we shall use the variables $(\beta,z)$ for points in $\d^2$ and the variables $(s,p)$ for points in $G$, so that
\[
s=\beta +\bar \beta z,\qquad p=z,\qquad \beta,z \in \d.
\]

Flat co-ordinates provide another characterization of domains biholomorphic to $G$.

The following lemma is a consequence of the Chain Rule.
\begin{lem}\label{lem10}
If $f=f(s,p)$ is a differentiable function on $G$, $\eta$ is defined on $\d^2$ as in equation \eqref{30} and $\xi =f \circ \eta$, then the following relations hold.
\be\label{40}
\der{\xi}{\beta} =\der{f}{s}
+\bar z \der{f}{\bar s},
\ee
\be\label{50}
\der{\xi}{\bar\beta} =z\der{f}{s}
+\der{f}{\bar s},
\ee
\be\label{60}
 \der{\xi}{z} =\bar\beta\der{f}{s}
+\der{f}{p},
\ee
\be\label{70}
\der{\xi}{\bar z} =\beta\der{f}{\bar s}
+\der{f}{\bar p}.
\ee
\end{lem}
\begin{thm}\label{thm30}
If $\Omega$ is a domain in $\c^2$, then $\Omega$ is biholomorphic to $G$ if and only if there exists a differentiable homeomorphism $\Xi=(\xi_1,\xi_2)$ from $\d^2$ onto $\Omega$ satisfying
\be\label{80}
\frac{\partial \xi_i}{\partial\bar\beta}=z\der{\xi_i}{\beta},\qquad i=1,2
\ee
and
\be\label{90}
\der{\xi_i}{\bar z} =0,\qquad i=1,2
\ee
at all $(\beta,z)\in \d^2$.
\end{thm}
\begin{proof}
First assume that $F\in \Omega(G)$ is a biholomorphic map of $G$ onto $\Omega$ and let $\Xi=F \circ \eta$. Since $\eta$ is a smooth homeomorphism of $\d^2$ onto $G$, $\Xi$ is a smooth homeomorphism of $\d^2$ onto $\Omega$.

If we set $F=(f_1,f_2)$ and $\Xi=(\xi_1,\xi_2)$, then $f_i$ is holomorphic and $\xi_i = f_i \circ \eta$ for $i=1,2$. Hence, using equations \eqref{40} and \eqref{50}, we see that
\[
\frac{\partial \xi_i}{\partial\bar\beta}=z\der{f_i}{s}=z\frac{\partial \xi_i}{\partial \beta},\qquad i=1,2,
\]
which proves that equation \eqref{80} holds. Also, equation \eqref{70} implies that the relation \eqref{90} holds.

Now assume that $\Xi=(\xi_1,\xi_2)$ is a differentiable homeomorphism  from $\d^2$ onto $\Omega$ satisfying equations \eqref{80} and \eqref{90}. Define $F=(f_1,f_2)$ by $F=\Xi \circ \eta^{-1}$. Since $\eta$ is a differentiable homeomorphism of $\d^2$ onto $G$, it follows that $F$ is a differentiable homeomorphism of $G$ onto $\Omega$. There remains to show that $F$ is holomorphic.

Since $\xi_i = f_i\circ \eta$,  we have
\begin{align*}
z\der{f_i}{s}+\der{f_i}{\bar s}&=   \frac{\partial \xi_i}{\partial\bar\beta} =  z\der{\xi_i}{\beta} \quad & \mbox{ by equations \eqref{50} and \eqref{80}} \\
	&=z\left(\der{f_i}{s}  +\bar z \der{f_i}{\bar s}\right) \quad& \mbox{ by equation \eqref{40}}.
\end{align*}
Thus
\[
(1-|z|^2)\der{f_i}{\bar s} = 0
\]
on $G$.  Since  $|z|<1$ when $(s,p) \in G$ it follows that 
\[
\frac{\partial f_i}{\partial \bar s}=0
\]
 throughout $G$. Hence $f_1,f_2$ are holomorphic on $G$.
\end{proof}

\section{Asymmetry of domains} \label{asymmetry}
\'E. Cartan's classification theorem \cite{ecartan}
\index{Cartan!\'E.}
 is based on his theory of {\em symmetric} spaces, 
\index{symmetric}
in the sense of the first paragraph of the paper.
\index{domain!bounded symmetric homogeneous}
In $\c^2$ and $\c^3$  (but not $\c^4$)  {\em every} bounded homogeneous domain is a symmetric space \cite{ecartan,fuks}.  In contrast, none of the `almost homogeneous' domains that we consider is symmetric. 

Let us say that a point $\la$ in a domain $\Omega$ is a {\em point of symmetry} of $\Omega$ if there exists a holomorphic self-map $\ga$ of $\Omega$ such that $\ga\circ\ga=\id_{\Omega}$ and $\la$ is an isolated fixed point of $\ga$.  Thus a domain is symmetric if every point of the domain is a point of symmetry.

From the fact that the automorphisms of the annulus  $\mathbb{A}_q$
\index{annulus}
\index{$\mathbb{A}_q$}
 are the maps $\omega z$ and $\omega z\inv$ for $\omega\in\t$ (for example, \cite[Theorem 6.2]{fisher}), it is easy to see that the only points of symmetry
in $\mathbb{A}_q$ are the points of the unit circle.
Hence $\mathbb{A}_q$ is not a symmetric domain.

\begin{proposition}
Neither the symmetrized bidisc nor the tetrablock contains a point of symmetry.
\index{tetrablock}
\end{proposition}
\begin{proof}
We sketch the proof for the tetrablock; that for the symmetrized bidisc is similar but simpler.

Let $E$ denote the tetrablock  defined in equation \eqref{deftet}.  
Every orbit in $E$ contains a point of the form $(0,0,p)$
 \cite[Theorem 5.2]{Y08}, so it suffices to show that no such point is a point of symmetry.  By \cite[Theorem 2.2]{awy}, the tetrablock is foliated by the `flat geodesics'
\[
C_{\beta_1\beta_2} \df \{(\beta_1+\bar\beta_2 z, \beta_2+\bar\beta_1z,z): z\in\d\}
\]
where $|\beta_1|+|\beta_2| < 1$.
These geodesics are permuted by the automorphisms of $E$ \cite[Theorem 5.1]{Y08}.  Moreover the `royal variety' $\{x\in E: x^1x^2=x^3\}$ is invariant under all automorphisms of $E$ (see the proof of \cite[Theorem 4.1]{Y08}).

Consider a holomorphic involution $\ga$ of $E$ that fixes $(0,0,p)$.  Then $\ga$ fixes the  flat geodesic containing $(0,0,p)$, which is $C_{00}$. Hence    $\ga$ fixes the only common point of $C_{00}$ and the royal variety, which is easily seen to be $(0,0,0)$.  It is shown in \cite[Proof of Theorem 4.1, foot of page 766]{Y08} that an automorphism $\ga$ of $E$  fixes $(0,0,0)$
if and only if either
\beq\label{case1}
\ga(x)=(\omega_1x^1, \omega_2 x^2, \omega_1\omega_2 x^3)
\eeq
or
\beq\label{case2}
\ga(x)=(\omega_2x^2, \omega_1 x^1, \omega_1\omega_2 x^3)
\eeq
 for some $\omega_1,\ \omega_2 \in\t$.

In the case that $\ga$ is of the form \eqref{case1}, since $\ga$ is an involution, we have $\omega_1^2=\omega_2^2 =1$.  Thus the four involutions of this form that fix $(0,0,p)$ have fixed points as in the following table.
\begin{center}
\begin{tabular}{cccc}
$\omega_1$ &$ \omega_2$ & $\ga(x)$ & Fixed points \\ \cline{1-4}
$1$ 		&$ 1$		& $ x$ 		& $E$ \\
$-1$ 	& $1$ 		& $(-x^1,x^2,-x^3)$ & $(0,z,0)$ \\
$1$ 		&$-1$		& $(x^1,-x^2,-x^3)$ & $(z,0,0)$\\
$-1$		&$-1$		&$(-x^1,-x^2,x^3)$  &$(0,0,z)$\\
\end{tabular}
\end{center}
where $z$ ranges over $\d$.  Hence in the case \eqref{case1}, $(0,0,p)$ is not an isolated fixed point of $\ga$.

In case \eqref{case2}, 
\[
\ga\circ\ga(x)=(\omega_1\omega_2x^1, \omega_1\omega_2x^2, (\omega_1\omega_2)^2x^3)
\]
and so the involutory property of $\ga$ corresponds to the condition $\omega_1\omega_2=1$.  Hence $\ga(x)=(\bar\omega x^2,\omega x^1,x^3)$
for some $\omega\in\t$.  Then  the fixed points of $\ga$ are the points $(x^1, \omega x^1, x^3)$ in $E$.   Hence
$(x^1,\omega x^1,p)$ is a fixed point of $\ga$ for all $x^1$ in a neighborhood of $0$, and so
 $(0,0,p)$ is not an isolated fixed point of $\ga$.
\end{proof}

\bibliography{references}

\printindex
\end{document}